\documentclass[final,3p]{elsarticle}
 \usepackage{graphics}
 \usepackage{graphicx}
 \usepackage{epsfig}
\usepackage{amssymb}
 \usepackage{amsthm}
 \usepackage{lineno}
 \usepackage{amsmath}
   \numberwithin{equation}{section}
\usepackage{mathrsfs}

\journal{Springer} 

\newtheorem{thm}{Theorem}[section]

\newtheorem{lem}[thm]{Lemma}

\biboptions{sort&compress,square}
\allowdisplaybreaks
\begin{document}
\begin{frontmatter}
\author{Tong Wu$^{a,b}$}
\ead{wut977@nenu.edu.cn}
\author{Yong Wang$^{c,*}$}
\ead{wangy581@nenu.edu.cn}
\cortext[cor]{Corresponding author.}
\address{$^a$Department of Mathematics, Northeastern University, Shenyang, 110819, China}
\address{$^b$Key Laboratory of Data Analytics and Optimization for Smart Industry\\ (Northeastern 
University), Ministry of Education, China}
\address{$^c$School of Mathematics and Statistics, Northeast Normal University,
Changchun, 130024, China}

\title{A general Dabrowski-Sitarz-Zalecki type theorems for manifold with boundary}
\begin{abstract}
In \cite{Wu}, we obtained the spectral Einstein functional associated with the Dirac operator for n-dimensional manifolds without boundary.  In this paper, we give the proof of general Dabrowski-Sitarz-Zalecki type theorems for the spectral Einstein functional associated with the Dirac operator on even and odd dimensional manifolds with boundary.
\end{abstract}
\end{frontmatter}
\section{Introduction}
 Until now, many geometers have studied noncommutative residues. In \cite{Gu,Wo}, authors found noncommutative residues are of great importance to the study of noncommutative geometry. In \cite{Co1}, Connes used the noncommutative residue to derive a conformal 4-dimensional Polyakov action analogy. Connes showed us that the noncommutative residue on a compact manifold $M$ coincided with the Dixmier's trace on pseudodifferential operators of order $-{\rm {dim}}M$ in \cite{Co2}.
And Connes claimed the noncommutative residue of the square of the inverse of the Dirac operator was proportioned to the Einstein-Hilbert action.  Kastler \cite{Ka} gave a
brute-force proof of this theorem. Kalau and Walze proved this theorem in the normal coordinates system simultaneously in \cite{KW} .
Ackermann proved that
the Wodzicki residue  of the square of the inverse of the Dirac operator ${\rm  Wres}(D^{-2})$ in turn is essentially the second coefficient
of the heat kernel expansion of $D^{2}$ in \cite{Ac}.

On the other hand, Wang generalized the Connes' results to the case of manifolds with boundary in \cite{Wa1,Wa2},
and proved the Kastler-Kalau-Walze type theorem for the Dirac operator and the signature operator on lower-dimensional manifolds
with boundary \cite{Wa3}. In \cite{Wa3,Wa4}, Wang computed $\widetilde{{\rm Wres}}[\pi^+D^{-1}\circ\pi^+D^{-1}]$ and $\widetilde{{\rm Wres}}[\pi^+D^{-2}\circ\pi^+D^{-2}]$, where the two operators are symmetric, in these cases the boundary term vanished. But for $\widetilde{{\rm Wres}}[\pi^+D^{-1}\circ\pi^+D^{-3}]$, J. Wang and Y. Wang got a nonvanishing boundary term \cite{Wa5}, and give a theoretical explanation for gravitational action on boundary. In others words, Wang provided a kind of method to study the Kastler-Kalau-Walze type theorem for manifolds with boundary. In \cite{DL}, the authors defined bilinear functionals of vector fields and differential forms, the densities of which yielded the  metric and Einstein tensors on even-dimensional Riemannian manifolds. In \cite{Wu}, the authors defined the spectral Einstein functional associated with the Dirac operator for manifolds with boundary, and compute the noncommutative residue $\widetilde{{\rm Wres}}[\pi^+(\nabla_X^{S(TM)}\nabla_Y^{S(TM)}D^{-2})\circ\pi^+(D^{-2})]$ on 4-dimensional compact manifolds. Motivated by \cite{Wa6,Wa7,Wu}, we compute the noncommutative residue $\widetilde{{\rm Wres}}[\pi^+(\nabla_X^{S(TM)}\nabla_Y^{S(TM)}D^{-2})\circ\pi^+(D^{-(n-2)})]$  and $\widetilde{{\rm Wres}}[\pi^+(\nabla_X^{S(TM)}\nabla_Y^{S(TM)}D^{-1})\circ\pi^+(D^{-(n-1)})]$ on even and odd dimensional compact manifolds. Our main theorems are as follows.
\\
\begin{thm}\label{thm1.11} If $M$ is an $n$-dimensional compact oriented spin manifolds without boundary, and $n$ is even, then we get the following equality:
\begin{align}
\label{1111}
{\rm Wres}[\nabla_X^{S(TM)}\nabla_Y^{S(TM)}D^{-n}]
&=\frac{(2\pi)^{\frac{n}{2}}}{3(\frac{n}{2}-1)!}\int_{M}[Ric(X,Y)-\frac{1}{2}sg(X,Y)]d{\rm Vol_{M}}+\frac{(2\pi)^{\frac{n}{2}}}{4(\frac{n}{2}-1)!}\int_{M}sg(X,Y)d{\rm Vol_{M}},
\end{align}
where $s$ is the scalar curvature, ${\rm Vol_{M}}$ is the volume of $M$ and $Ric$ denotes Ricci tensor on $M$.
\end{thm}
\begin{thm}\label{thm12222}
Let $M$ be an $n$-dimensional oriented
compact spin manifold with boundary $\partial M$ and n is even, then we get the following equality:
\begin{align}
\label{122222}
&\widetilde{{\rm Wres}}[\pi^+(\nabla_X^{S(TM)}\nabla_Y^{S(TM)}D^{-2})\circ\pi^+(D^{-(n-2)})]\nonumber\\
&=\frac{(2\pi)^{\frac{n}{2}}}{3(\frac{n}{2}-1)!}\int_{M}[Ric(X,Y)-\frac{1}{2}sg(X,Y)]d{\rm Vol_{M}}+\frac{(2\pi)^{\frac{n}{2}}}{4(\frac{n}{2}-1)!}\int_{M}sg(X,Y)d{\rm Vol_{M}}-\int_{\partial M}\frac{1}{2}Vol(S^{n-2})h'(0)2^{\frac{n}{2}}\nonumber\\
&\bigg\{(\frac{n}{2}-1)\bigg(\frac{1}{n-1}g(X^T,Y^T)+X_nY_n\bigg)\frac{\pi}{2(\frac{n}{2}+2)!}A_{0}+(\frac{n}{2}-1)X_nY_n\frac{\pi i}{(\frac{n}{2}+2)!}A_{1}+(\frac{n}{2}-1)\bigg[\frac{i}{n-1}\bigg(h'(0)\nonumber\\
&g(X^T,Y^T)+\partial_{x_n}(g(X^T,Y^T))\bigg)-\partial_{x_n}(X_nY_n)\bigg]\frac{\pi i}{(\frac{n}{2}+1)!}A_2-(1-\frac{n}{2})\bigg(\frac{1}{n-1}g(X^T,Y^T)-X_nY_n\bigg)\nonumber\\
&\frac{2\pi i}{(\frac{n}{2}+2)!}B_0-\bigg(\frac{i}{n-1}g(X^T,Y^T)+X_nY_n\bigg)\frac{\pi i}{2(\frac{n}{2}+2)!}C_{0}-\frac{1}{n-1}g(X^T,Y^T)(1-\frac{n}{2})\frac{2\pi }{(1+\frac{n}{2})!}D_0\nonumber\\
&+\frac{1}{n-1}g(X^T,Y^T)(1-\frac{n}{2})\frac{4\pi }{(2+\frac{n}{2})!}D_1-X_nY_n(1-\frac{n}{2})\frac{32\pi i }{(\frac{n}{2}+1)!}D_2\bigg\}d{\rm Vol_{M}},
\end{align}
where $A_0$-$D_2$ are defined in section 3.
\end{thm}
\begin{thm}\label{mmmmm}
Let $M$ be an $n-1$ dimensional
oriented compact spin manifold with boundary $\partial M$ and $n-1$ is odd, then
\begin{align}\label{8888}
&\widetilde{{\rm Wres}}[\pi^+(\nabla_X^{S(TM)}\nabla_Y^{S(TM)}D^{-2})\circ\pi^+(D^{-(n-2)})]\nonumber\\
&=\int_{\partial M}-Vol(S^{n-2})\bigg(\frac{1}{n-1}g(X^T,Y^T)+X_nY_n\bigg)(1-\frac{n}{2})2^{\frac{n}{2}-1}\frac{2\pi i}{\frac{n}{2}!}M_0d{\rm Vol_{M}},\nonumber\\
\end{align}
where $M_0$ is defined in section 4.\\
\end{thm}
\indent The paper is organized in the following way. In Section \ref{section:2}, we define the spectral Einstein functional associated with the Dirac operator and and get the noncommutative residue for manifolds without boundary. In Section \ref{section:3}, we prove the Dabrowski-Sitarz-Zalecki type theorem for the spectral Einstein functional associated with the Dirac operator on even dimensional manifolds with boundary. In Section \ref{section:4}, we prove the Dabrowski-Sitarz-Zalecki type theorem for the spectral Einstein functional associated with the Dirac operator on odd dimensional manifolds with boundary.
\section{The spectral Einstein functional associated with the Dirac operator}
\label{section:2}
Firstly we recall the definition of Dirac operator. Let $M$ be an $n$-dimensional ($n\geq 4$) even oriented compact Riemannian manifold with a Riemannian metric $g^{M}$ and let $\nabla^L$ be the Levi-Civita connection about $g^{M}$. In the fixed orthonormal frame $\{e_1,\cdots,e_n\}$, the connection matrix $(\omega_{s,t})$ is defined by
\begin{equation}
\label{a2}
\nabla^L(e_1,\cdots,e_n)= (e_1,\cdots,e_n)(\omega_{s,t}).
\end{equation}
\indent Let $\epsilon (e_j^*)$,~$\iota (e_j^*)$ be the exterior and interior multiplications respectively, where $e_j^*=g^{TM}(e_j,\cdot)$.
Write
\begin{equation}
\label{a3}
\widehat{c}(e_j)=\epsilon (e_j^* )+\iota
(e_j^*);~~
c(e_j)=\epsilon (e_j^* )-\iota (e_j^* ),
\end{equation}
which satisfies
\begin{align}
\label{a4}
&\widehat{c}(e_i)\widehat{c}(e_j)+\widehat{c}(e_j)\widehat{c}(e_i)=2g^{M}(e_i,e_j);~~\nonumber\\
&c(e_i)c(e_j)+c(e_j)c(e_i)=-2g^{M}(e_i,e_j);~~\nonumber\\
&c(e_i)\widehat{c}(e_j)+\widehat{c}(e_j)c(e_i)=0.\nonumber\\
\end{align}
By \cite{Y}, we have the Dirac operator
\begin{align}
\label{a5}
&D=\sum^n_{i=1}c(e_i)[e_i-\frac{1}{4}\sum_{s,t}\omega_{s,t}
(e_i)c(e_s)c(e_t)].\nonumber\\
\end{align}
\indent We define $\nabla_X^{S(TM)}:=X+\frac{1}{4}\Sigma_{ij}\langle\nabla_X^L{e_i},e_j\rangle c(e_i)c(e_j)$, which is a spin connection. Set\\ $A(X)=\frac{1}{4}\Sigma_{ij}\langle\nabla_X^L{e_i},e_j\rangle c(e_i)c(e_j)$, then
\begin{align}\label{ddd}
\nabla_X^{S(TM)}\nabla_Y^{S(TM)}&=[X+A(X)][Y+A(Y)]\nonumber\\
&=XY+X\cdot A(Y)+A(X)Y+A(X)A(Y)\nonumber\\
&=XY+X[A(Y)]+A(Y)X+A(X)Y+A(X)A(Y),\nonumber\\
\end{align}
where $X=\sum_{j=1}^nX_j\partial_{x_j}, Y=\sum_{l=1}^nY_l\partial_{x_l}$ and $XY=\sum_{j,l=1}^n(X_jY_l\partial_{x_j}\partial_{x_l}+X_j\frac{\partial{Y_l}}{\partial_{x_j}}\partial_{x_l}).$ \\
\indent Let $g^{ij}=g(dx_{i},dx_{j})$, $\xi=\sum_{k}\xi_{j}dx_{j}$ and $\nabla^L_{\partial_{i}}\partial_{j}=\sum_{k}\Gamma_{ij}^{k}\partial_{k}$,  we denote that
\begin{align}
&\sigma_{i}=-\frac{1}{4}\sum_{s,t}\omega_{s,t}
(e_i)c(e_i)c(e_s)c(e_t)
;~~~\xi^{j}=g^{ij}\xi_{i};~~~~\Gamma^{k}=g^{ij}\Gamma_{ij}^{k};~~~~\sigma^{j}=g^{ij}\sigma_{i};
~~~~a^{j}=g^{ij}a_{i}.
\end{align}
And by $\partial_{x_j}=-\sqrt{-1}\xi_j$, we have the following lemmas.
\begin{lem}\label{lem3} The following identities hold:
\begin{align}
\label{b22}
&\sigma_{0}(\nabla_X^{S(TM)}\nabla_Y^{S(TM)})=X[A(Y)]+A(X)A(Y);\nonumber\\
&\sigma_{1}(\nabla_X^{S(TM)}\nabla_Y^{S(TM)})=\sqrt{-1}\sum_{j,l=1}^nX_j\frac{\partial{Y_l}}{\partial_{x_j}}\xi_l+\sqrt{-1}\sum_jA(Y)X_j\xi_j+\sqrt{-1}\sum_lA(Y)Y_l\xi_l;\nonumber\\
&\sigma_{2}(\nabla_X^{S(TM)}\nabla_Y^{S(TM)})=-\sum_{j,l=1}^nX_jY_l\xi_j\xi_l.\nonumber\\
\end{align}
\end{lem}
\begin{lem}\cite{Ka}\label{lem356} The following identities hold:
\begin{align}
\label{b22222}
&\sigma_{-2}(D^{-2})=|\xi|^{-2};\nonumber\\
&\sigma_{-3}(D^{-2})=-\sqrt{-1}|\xi|^{-4}\xi_k(\Gamma^k-2\delta^k)-\sqrt{-1}|\xi|^{-6}2\xi^j\xi_\alpha\xi_\beta\partial_jg^{\alpha\beta}.\nonumber\\
\end{align}
\end{lem}
\begin{lem}\cite{Wa5}\label{lem9356} The following identities hold:
\begin{align}
\label{bii}
&\sigma_{-1}(D^{-1})=\frac{\sqrt{-1}c(\xi)}{|\xi|^{2}};\nonumber\\
&\sigma_{-2}(D^{-1})=\frac{c(\xi)\sigma_0(D)c(\xi)}{|\xi|^4}
+\frac{c(\xi)}{|\xi|^6}\Sigma_jc(\mathrm{d}x_j)(\partial_{x_j}[c(\xi)]|\xi|^2-c(\xi)\partial_{x_j}(|\xi|^2)).\nonumber\\
\end{align}
\end{lem}
\indent Write
 \begin{eqnarray}
D_x^{\alpha}&=(-i)^{|\alpha|}\partial_x^{\alpha};
~\sigma(D_t)=p_1+p_0;
~(\sigma(D_t)^{-1})=\sum^{\infty}_{j=1}q_{-j}.
\end{eqnarray}

\indent By the composition formula of pseudodifferential operators, we have
\begin{align}
1=\sigma(D\circ D^{-1})&=\sum_{\alpha}\frac{1}{\alpha!}\partial^{\alpha}_{\xi}[\sigma(D)]
D_x^{\alpha}[\sigma(D^{-1})]\nonumber\\
&=(p_1+p_0)(q_{-1}+q_{-2}+q_{-3}+\cdots)\nonumber\\
&~~~+\sum_j(\partial_{\xi_j}p_1+\partial_{\xi_j}p_0)(
D_{x_j}q_{-1}+D_{x_j}q_{-2}+D_{x_j}q_{-3}+\cdots)\nonumber\\
&=p_1q_{-1}+(p_1q_{-2}+p_0q_{-1}+\sum_j\partial_{\xi_j}p_1D_{x_j}q_{-1})+\cdots,
\end{align}
so
\begin{equation}
q_{-1}=p_1^{-1};~q_{-2}=-p_1^{-1}[p_0p_1^{-1}+\sum_j\partial_{\xi_j}p_1D_{x_j}(p_1^{-1})].
\end{equation}
Then, we have the following lemma
\begin{lem}\label{lema}The following identities hold:
\begin{align}\label{mki}
\sigma_{0}(\nabla_X^{S(TM)}\nabla_Y^{S(TM)}D^{-2})&=-\sum_{j,l=1}^nX_jY_l\xi_j\xi_l|\xi|^{-2};\nonumber\\
\sigma_{1}(\nabla_X^{S(TM)}\nabla_Y^{S(TM)}D^{-1})&=-i\sum_{j,l=1}^nX_jY_l\xi_j\xi_lc(\xi)|\xi|^{-2}.
\end{align}
\end{lem}
By \cite{DL}, we have the following theorem
\begin{thm}\label{thm2} If $M$ is an $n$-dimensional compact oriented manifolds without boundary, and $n$ is even, then we get the following equality:
\begin{align}
\label{a29}
{\rm Wres}[\nabla_X^{S(TM)}\nabla_Y^{S(TM)}D^{-n}]
&=\frac{(2\pi)^{\frac{n}{2}}}{3(\frac{n}{2}-1)!}\int_{M}[Ric(X,Y)-\frac{1}{2}sg(X,Y)]d{\rm Vol_{M}}+\frac{(2\pi)^{\frac{n}{2}}}{4(\frac{n}{2}-1)!}\int_{M}sg(X,Y)d{\rm Vol_{M}},
\end{align}
where $s$ is the scalar curvature, ${\rm Vol_{M}}$ is the volume of $M$ and $Ric$ denotes Ricci tensor on $M$.
\end{thm}

\section{ Dabrowski-Sitarz-Zalecki type theorems for even dimensional manifolds with boundary}
\label{section:3}
 In this section, we compute the noncommutative residue $\widetilde{{\rm Wres}}[\pi^+(\nabla_X^{S(TM)}\nabla_Y^{S(TM)}D^{-2})\circ\pi^+(D^{-(n-2)})]$ and $\widetilde{{\rm Wres}}[\pi^+(\nabla_X^{S(TM)}\nabla_Y^{S(TM)}D^{-1})\circ\pi^+(D^{-(n-1)})]$ on $n$ dimensional oriented compact manifolds with boundary, and $n$ is even. Firstly, we recall some basic facts and formulas about Boutet de
Monvel's calculus and the definition of the noncommutative residue for manifolds with boundary which will be used in the following. For more details, see Section 2 in \cite{Wa3}.\\
 \indent Let $M$ be an n-dimensional compact oriented manifold with boundary $\partial M$, and $n$ is even.
We assume that the metric $g^{M}$ on $M$ has the following form near the boundary,
\begin{equation}
\label{b1}
g^{M}=\frac{1}{h(x_{n})}g^{\partial M}+dx _{n}^{2},
\end{equation}
where $g^{\partial M}$ is the metric on $\partial M$ and $h(x_n)\in C^{\infty}([0, 1)):=\{\widehat{h}|_{[0,1)}|\widehat{h}\in C^{\infty}((-\varepsilon,1))\}$ for
some $\varepsilon>0$ and $h(x_n)$ satisfies $h(x_n)>0$, $h(0)=1,$ where $x_n$ denotes the normal directional coordinate.\\
\indent Firstly, we compute $\widetilde{{\rm Wres}}[\pi^+(\nabla_X^{S(TM)}\nabla_Y^{S(TM)}D^{-2})\circ\pi^+(D^{-(n-2)})]$, similar to \cite{Wa3},  we can get
\begin{align}
\label{b14}
&\widetilde{{\rm Wres}}[\pi^+(\nabla_X^{S(TM)}\nabla_Y^{S(TM)}D^{-2})\circ\pi^+(D^{-(n-2)})]\nonumber\\
&=\int_M\int_{|\xi|=1}{\rm
trace}_{\wedge^*T^*M\bigotimes\mathbb{C}}[\sigma_{-n}(\nabla_X^{S(TM)}\nabla_Y^{S(TM)}D^{-2}\circ D^{-(n-2)})]\sigma(\xi)dx+\int_{\partial M}\Phi,
\end{align}
where
\begin{align}
\label{qqq15}
\Phi &=\int_{|\xi'|=1}\int^{+\infty}_{-\infty}\sum^{\infty}_{j, k=0}\sum\frac{(-i)^{|\alpha|+j+k+1}}{\alpha!(j+k+1)!}
\times {\rm trace}_{\wedge^*T^*M\bigotimes\mathbb{C}}[\partial^j_{x_n}\partial^\alpha_{\xi'}\partial^k_{\xi_n}\sigma^+_{r}(\nabla_X^{S(TM)}\nabla_Y^{S(TM)}D^{-2})(x',0,\xi',\xi_n)
\nonumber\\
&\times\partial^\alpha_{x'}\partial^{j+1}_{\xi_n}\partial^k_{x_n}\sigma_{l}(D^{-(n-2)})(x',0,\xi',\xi_n)]d\xi_n\sigma(\xi')dx',
\end{align}
and the sum is taken over $r+l-k-j-|\alpha|-1=-n,~~r\leq 0,~~l\leq -(n-2)$.\\

\indent The interior of $\widetilde{{\rm Wres}}[\pi^+(\nabla_X^{S(TM)}\nabla_Y^{S(TM)}D^{-2})\circ\pi^+(D^{-(n-2)})]$ has been given in Theorem \ref{thm2}. Therefore, we only need to compute $\int_{\partial M} \Phi$. When $n$ is even, then ${\rm tr}_{S(TM)}[{\rm \texttt{id}}]={\rm dim}(\wedge^*(\mathbb{R}^2))=2^{\frac{n}{2}}$, the sum is taken over $
r+l-k-j-|\alpha|=-n+1,~~r\leq 0,~~l\leq -(n-2),$ then we have the following five cases:
~\\
\noindent  {\bf case a)~I)}~$r=0,~l=-(n-2),~k=j=0,~|\alpha|=1$.\\
\noindent By (\ref{qqq15}), we get
\begin{equation}
\label{b24}
\Phi_1=-\int_{|\xi'|=1}\int^{+\infty}_{-\infty}\sum_{|\alpha|=1}
 {\rm trace}[\partial^\alpha_{\xi'}\pi^+_{\xi_n}\sigma_{0}(\nabla_X^{S(TM)}\nabla_Y^{S(TM)}D^{-2})\times
 \partial^\alpha_{x'}\partial_{\xi_n}\sigma_{-(n-2)}(D^{-(n-2)})](x_0)d\xi_n\sigma(\xi')dx'.
\end{equation}
By Lemma 2.2 in \cite{Wa3}, for $i<n$, then
\begin{equation}
\label{b25}
\partial_{x_i}\sigma_{-(n-2)}(D^{-(n-2)})(x_0)=\partial_{x_i}{(|\xi|^{(2-n)})}(x_0)
=\partial_{x_i}(|\xi|^2)^{(1-\frac{n}{2})}(x_0)=(1-\frac{n}{2}) (|\xi|^2)^{-\frac{n}{2}}
 \partial_{x_i}(|\xi|^2)(x_0)=0,
\end{equation}
\noindent so $\Phi_1=0$.\\
 \noindent  {\bf case a)~II)}~$r=0,~l=-(n-2),~k=|\alpha|=0,~j=1$.\\
\noindent By (\ref{qqq15}), we get
\begin{equation}
\label{b26}
\Phi_2=-\frac{1}{2}\int_{|\xi'|=1}\int^{+\infty}_{-\infty} {\rm
trace} [\partial_{x_n}\pi^+_{\xi_n}\sigma_{0}(\nabla_X^{S(TM)}\nabla_Y^{S(TM)}D^{-2})\times
\partial_{\xi_n}^2\sigma_{-(n-2)}(D^{-(n-2)})](x_0)d\xi_n\sigma(\xi')dx'.
\end{equation}
\noindent By (3.16) in \cite{Wa6}, we have\\
\begin{align}\label{b237}
\partial^2_{\xi_n} \sigma_{-(n-2)}(D^{-(n-2)})(x_0)
&=\partial^2_{\xi_n}\big((|\xi|^2)^{1-\frac{n}{2}}\big) (x_0)
=\partial_{\xi_n}\Big( (1-\frac{n}{2}) (|\xi|^{2})^{-\frac{n}{2}}  \partial_{\xi_n} (|\xi|^2)\Big) (x_0)\nonumber\\
&=(1-\frac{n}{2})(-\frac{n}{2})(|\xi|^{2})^{-\frac{n}{2}-1} \big(\partial_{\xi_n}|\xi|^2\big)^{2}(x_0)+(1-\frac{n}{2})(|\xi|^{2})^{-\frac{n}{2}} \partial^2_{\xi_n}(|\xi|^2 (x_0))\nonumber\\
&=\Big((2n-2) \xi_n^{2}-2 \Big)(\frac{n}{2}-1)(1+\xi_n^{2})^{(-\frac{n}{2}-1)}.
\end{align}
By Lemma \ref{lem3}, we have\\
\begin{align}\label{b27}
\partial_{x_n}\sigma_{0}(\nabla_X^{S(TM)}\nabla_Y^{S(TM)}D^{-2})&=\partial_{x_n}\bigg(-\sum_{j,l=1}^nX_jY_l\xi_j\xi_l|\xi|^{-2}\bigg)\nonumber\\
&=\frac{\sum_{j,l=1}^nX_jY_l\xi_j\xi_lh'(0)|\xi'|^2}{(1+\xi_n^2)^2}-\frac{\sum_{j,l=1}^n\frac{\partial{X_j}}{\partial_{x_n}}Y_l\xi_j\xi_l}{1+\xi_n^2}-\frac{\sum_{j,l=1}^nX_j\frac{\partial{Y_l}}{\partial_{x_n}}\xi_j\xi_l}{1+\xi_n^2}.
\end{align}
Then, we have
\begin{align}\label{b28}
&\partial_{x_n}\pi^+_{\xi_n}\sigma_{0}(\nabla_X^{S(TM)}\nabla_Y^{S(TM)}D^{-2})\nonumber\\
&=\pi^+_{\xi_n}\partial_{x_n}\sigma_{0}(\nabla_X^{S(TM)}\nabla_Y^{S(TM)}D^{-2})\nonumber\\
&=-\frac{i\xi_n}{4(\xi_n-i)^2}\sum_{j,l=1}^{n-1}X_jY_l\xi_j\xi_lh'(0)+\frac{2-i\xi_n}{4(\xi_n-i)^2}X_nY_nh'(0)-\frac{i}{4(\xi_n-i)^2}\sum_{j=1}^{n-1}\bigg(X_jY_n+X_nY_j\bigg)\xi_j\nonumber\\
&+\frac{i}{2(\xi_n-i)}\sum_{j,l=1}^{n-1}\bigg(\frac{\partial{X_j}}{\partial_{x_n}}Y_l\xi_j\xi_l+X_j\frac{\partial{Y_l}}{\partial_{x_n}}\xi_j\xi_l\bigg)-\frac{i}{2(\xi_n-i)}\bigg(X_n\frac{\partial{Y_n}}{\partial_{x_n}}+\frac{\partial{X_n}}{\partial_{x_n}}Y_n\bigg)\nonumber\\
&-\frac{1}{2(\xi_n-i)}\sum_{j=1}^{n-1}\bigg(\frac{\partial{X_j}}{\partial_{x_n}}Y_n+\frac{\partial{X_n}}{\partial_{x_n}}Y_j+X_j\frac{\partial{Y_n}}{\partial_{x_n}}+X_n\frac{\partial{Y_j}}{\partial_{x_n}}\bigg)\xi_j.\nonumber\\
\end{align}
We note that $i<n,~\int_{|\xi'|=1}\xi_{i_{1}}\xi_{i_{2}}\cdots\xi_{i_{2d+1}}\sigma(\xi')=0$,
so we omit some items that have no contribution for computing {\bf case a)~II)}.\\
Let $X=X^T+X_n\partial_n,~Y=Y^T+Y_n\partial_n,$ then we have $\sum_{j=1}^{n-1}X_jY_j=g(X^T,Y^T).$ And by $\int_{S^{n-2}}\xi_j\xi_l\sigma(\xi')=\frac{1}{n-1}\delta_{jl}Vol(S^{n-2})$,  where $Vol(S^{n-2})$ is the canonical volume of $S^{n-2}$,
we have
\begin{align}\label{35}
&\int_{|\xi'|}\int_{-\infty}^{\infty}{\rm tr}\bigg[-i\sum_{j,l=1}^{n-1}X_jY_l\xi_j\xi_lh'(0)\frac{\xi_n}{4(\xi_n-i)^2}\times
 \Big((2n-2) \xi_n^{2}-2 \Big)(\frac{n}{2}-1)(1+\xi_n^{2})^{(-\frac{n}{2}-1)}\bigg]d\xi_n\sigma(\xi')dx'\nonumber\\
&=\int_{|\xi'|}-\frac{ih'(0)}{4}\sum_{j,l=1}^{n-1}X_jY_l\xi_j\xi_l(\frac{n}{2}-1)2^{\frac{n}{2}}   \int_{\Gamma^+}\frac{(2n-2)\xi_n^3-2\xi_n}
   {(\xi_n-i)^{(\frac{n}{2}+3)}(\xi_n+i)^{(\frac{n}{2}+1)}}d\xi_ndx'\nonumber\\
&=\int_{|\xi'|}-\frac{ih'(0)}{4}\sum_{j,l=1}^{n-1}X_jY_l\xi_j\xi_l(\frac{n}{2}-1)2^{\frac{n}{2}}\sigma(\xi') \frac{2\pi i}{(\frac{n}{2}+2)!}
\left[\frac{(2n-2)\xi_n^3-2\xi_n}{(\xi_n+i)^{(\frac{n}{2}+1)}}\right]
 ^{(\frac{n}{2}+2)}\bigg|_{\xi_n=i}dx'\nonumber\\
 &=Vol(S^{n-2})h'(0)g(X^T,Y^T)\frac{\frac{n}{2}-1}{n-1}2^{\frac{n}{2}}\frac{\pi}{2(\frac{n}{2}+2)!}\left[\frac{(2n-2)\xi_n^3-2\xi_n}{(\xi_n+i)^{(\frac{n}{2}+1)}}\right]
 ^{(\frac{n}{2}+2)}\bigg|_{\xi_n=i}dx'\nonumber\\
  &:=Vol(S^{n-2})h'(0)g(X^T,Y^T)\frac{\frac{n}{2}-1}{n-1}2^{\frac{n}{2}}\frac{\pi}{2(\frac{n}{2}+2)!}A_{0}dx'.
\end{align}
Similarly, we can obtain the integral of the other terms. Then, we get\\
\begin{align}\label{35kkk}
\Phi_2
&=-\frac{1}{2}Vol(S^{n-2})h'(0)2^{\frac{n}{2}}(\frac{n}{2}-1)\bigg\{-\frac{i}{2}\bigg(\frac{1}{n-1}g(X^T,Y^T)+X_nY_n\bigg)\frac{\pi i}{(\frac{n}{2}+2)!}A_{0}+X_nY_n\frac{\pi i}{(\frac{n}{2}+2)!}A_{1}\nonumber\\
&+\bigg[\frac{i}{n-1}\bigg(h'(0)g(X^T,Y^T)+\partial_{x_n}(g(X^T,Y^T))\bigg)-\partial_{x_n}(X_nY_n)\bigg]\frac{\pi i}{(\frac{n}{2}+1)!}A_2\bigg\}dx',
\end{align}
where let $C_N^K=\frac{N!}{K!(N-K)!}$, we have\\
\begin{align}
A_0&=\left[\frac{(2n-2)\xi_n^3-2\xi_n}{(\xi_n+i)^{(\frac{n}{2}+1)}}\right]
 ^{(\frac{n}{2}+2)}\bigg|_{\xi_n=i}\nonumber\\
 &=-i^{-n-4}2^{-n-3}\bigg((4n-4)C_{-\frac{n}{2}-3}^{\frac{n}{2}-1}+(6n-6)C_{-\frac{n}{2}-3}^{\frac{n}{2}}+(3n-2)C_{-\frac{n}{2}-3}^{\frac{n}{2}+1}+\frac{n}{2}C_{-\frac{n}{2}-3}^{\frac{n}{2}+2}\bigg)(\frac{n}{2}+2)!;\nonumber\\
A_1&=\left[\frac{(2n-2)\xi_n^2-2}{(\xi_n+i)^{(\frac{n}{2}+1)}}\right]
 ^{(\frac{n}{2}+2)}\bigg|_{\xi_n=i}\nonumber\\
  &=(2i)^{-n-3}\bigg((2n-2)C_{-\frac{n}{2}-3}^{\frac{n}{2}}+(2n-2)C_{-\frac{n}{2}-3}^{\frac{n}{2}+1}+\frac{n}{2}C_{-\frac{n}{2}-3}^{\frac{n}{2}+2}\bigg)(\frac{n}{2}+2)!;\nonumber\\
A_2&=\left[\frac{(2n-2)\xi_n^3-2\xi_n}{(\xi_n+i)^{(\frac{n}{2}+1)}}\right]
 ^{(\frac{n}{2}+1)}\bigg|_{\xi_n=i}\nonumber\\
  &=i^{-n-3}2^{-n}\bigg((4n-4)C_{-\frac{n}{2}-1}^{\frac{n}{2}-2}+(6n-6)C_{-\frac{n}{2}-1}^{\frac{n}{2}-1}+(3n-2)C_{-\frac{n}{2}-1}^{\frac{n}{2}}+\frac{n}{2}C_{-\frac{n}{2}-1}^{\frac{n}{2}+1}\bigg)(\frac{n}{2}+2)!.\nonumber\\
\end{align}
\noindent  {\bf case a)~III)}~$r=0,~l=-(n-2),~j=|\alpha|=0,~k=1$.\\
\noindent By (\ref{qqq15}), we get
\begin{align}\label{36}
\Phi_3&=-\frac{1}{2}\int_{|\xi'|=1}\int^{+\infty}_{-\infty}
{\rm trace} [\partial_{\xi_n}\pi^+_{\xi_n}\sigma_{0}(\nabla_X^{S(TM)}\nabla_Y^{S(TM)}D^{-2})\times
\partial_{\xi_n}\partial_{x_n}\sigma_{-(n-2)}(D^{-(n-2)})](x_0)d\xi_n\sigma(\xi')dx'\nonumber\\
&=\frac{1}{2}\int_{|\xi'|=1}\int^{+\infty}_{-\infty}
{\rm trace} [\partial_{\xi_n}^2\pi^+_{\xi_n}\sigma_{0}(\nabla_X^{S(TM)}\nabla_Y^{S(TM)}D^{-2})\times
\partial_{x_n}\sigma_{-(n-2)}(D^{-(n-2)})](x_0)d\xi_n\sigma(\xi')dx'.
\end{align}
By Lemma \ref{lem356}, we have
\begin{align}\label{aa38}
\pi^+_{\xi_n}\sigma_{0}(\nabla_X^{S(TM)}\nabla_Y^{S(TM)}D^{-2})=\frac{i}{2(\xi_n-i)}\sum_{j,l=1}^{n-1}X_jY_l\xi_j\xi_l-\frac{1}{2(\xi_n-i)}X_nY_n-\frac{1}{2(\xi_n-i)}\sum_{j=1}^{n-1}(X_jY_n+X_nY_j)\xi_j.
\end{align}
By (3.21) in \cite{Wa6}, we have 
\begin{eqnarray}\label{37}
\partial_{x_n} \big(\sigma_{-(n-2)}(D^{-(n-2)})\big)(x_0)
=\partial_{x_n}\big((|\xi|^2)^{1-\frac{n}{2}}\big) (x_0)
=h'(0)(1-\frac{n}{2})(1+\xi_{n}^{2})^{-\frac{n}{2}}.
\end{eqnarray}
Then by (\ref{aa38}), we get
\begin{align}\label{mmmmm}
&\partial_{\xi_n}^2\pi^+_{\xi_n}\sigma_{0}(\nabla_X^{S(TM)}\nabla_Y^{S(TM)}D^{-2})\nonumber\\
&=\frac{i}{(\xi_n-i)^3}\sum_{j,l=1}^{n-1}X_jY_l\xi_j\xi_l-\frac{i}{(\xi_n-i)^3}X_nY_n-\frac{i}{(\xi_n-i)^3}\sum_{j=1}^{n-1}(X_jY_n+X_nY_j)\xi_j.
\end{align}
Next, we perform the corresponding integral calculation on the above results.
\begin{align}\label{3oooo5}
&\int_{|\xi'|}\int_{-\infty}^{\infty}{\rm tr }\bigg[i\sum_{j,l=1}^{n-1}X_jY_l\xi_j\xi_l\frac{1}{(\xi_n-i)^3}\times
 h'(0)(1-\frac{n}{2})(1+\xi_n^{2})^{-\frac{n}{2}}\bigg]d\xi_n\sigma(\xi')dx'\nonumber\\
&= \int_{|\xi'|}ih'(0)\sum_{j,l=1}^{n-1}X_jY_l\xi_j\xi_l(1-\frac{n}{2})2^{\frac{n}{2}}\sigma(\xi')\int_{\Gamma^+}\frac{1}
   {(\xi_n-i)^{(\frac{n}{2}+3)}(\xi_n+i)^{(\frac{n}{2})}}d\xi_ndx'\nonumber\\
&= \int_{|\xi'|}ih'(0)\sum_{j,l=1}^{n-1}X_jY_l\xi_j\xi_l(1-\frac{n}{2})2^{\frac{n}{2}}\sigma(\xi') \frac{2\pi i}{(\frac{n}{2}+2)!}
\left[\frac{1}{(\xi_n+i)^{\frac{n}{2}}}\right]
 ^{(\frac{n}{2}+2)}\bigg|_{\xi_n=i}dx'\nonumber\\
 &=-Vol(S^{n-2})h'(0)g(X^T,Y^T)\frac{1-\frac{n}{2}}{n-1}2^{\frac{n}{2}}\frac{2\pi}{(\frac{n}{2}+2)!}\left[\frac{1}{(\xi_n+i)^{(\frac{n}{2})}}\right]
 ^{(\frac{n}{2}+2)}\bigg|_{\xi_n=i}dx'\nonumber\\
  &:=-Vol(S^{n-2})h'(0)g(X^T,Y^T)\frac{1-\frac{n}{2}}{n-1}2^{\frac{n}{2}}\frac{2\pi}{(\frac{n}{2}+2)!}B_{0}dx',
\end{align}
and
\begin{align}\label{311115}
& \int_{|\xi'|}\int_{-\infty}^{\infty}{\rm tr }\bigg[-X_nY_n\frac{1}{(\xi_n-i)^3}\times
 h'(0)(1-\frac{n}{2})(1+\xi_n^{2})^{-\frac{n}{2}}\bigg]d\xi_n\sigma(\xi')dx'\nonumber\\
&= \int_{|\xi'|}-h'(0)X_nY_n(1-\frac{n}{2})2^{\frac{n}{2}}\sigma(\xi')\int_{\Gamma^+}\frac{1}
   {(\xi_n-i)^{(\frac{n}{2}+3)}(\xi_n+i)^{(\frac{n}{2})}}d\xi_ndx'\nonumber\\
&= \int_{|\xi'|}-h'(0)X_nY_n(1-\frac{n}{2})2^{\frac{n}{2}}\sigma(\xi') \frac{2\pi i}{(\frac{n}{2}+2)!}
\left[\frac{1}{(\xi_n+i)^{\frac{n}{2}}}\right]
 ^{(\frac{n}{2}+2)}\bigg|_{\xi_n=i}dx'\nonumber\\
 &=-Vol(S^{n-2})h'(0)X_nY_n(1-\frac{n}{2})2^{\frac{n}{2}}\frac{2\pi i}{(\frac{n}{2}+2)!}\left[\frac{1}{(\xi_n+i)^{(\frac{n}{2})}}\right]
 ^{(\frac{n}{2}+2)}\bigg|_{\xi_n=i}dx'\nonumber\\
  &:=-Vol(S^{n-2})h'(0)X_nY_n(1-\frac{n}{2})2^{\frac{n}{2}}\frac{2\pi i}{(\frac{n}{2}+2)!}B_{0}dx'.
\end{align}
Therefore, we get
\begin{align}\label{41}
\Phi_3&=Vol(S^{n-2})h'(0)(1-\frac{n}{2})2^{\frac{n}{2}}\bigg(\frac{1}{n-1}g(X^T,Y^T)-X_nY_n\bigg)\frac{\pi i}{(\frac{n}{2}+2)!}B_0dx',
\end{align}
where let $A_N^K=\frac{N!}{(N-K)!}$, we have\\
\begin{align}
B_0&=\left[\frac{1}{(\xi_n+i)^{\frac{n}{2}}}\right]
 ^{(\frac{n}{2}+2)}\bigg|_{\xi_n=i}=(2i)^{-n-2}A_{-\frac{n}{2}}^{\frac{n}{2}+2}.\nonumber\\
\end{align}
\noindent  {\bf case b)}~$r=0,~l=-(n-1),~k=j=|\alpha|=0$.\\
\noindent By (\ref{qqq15}), we get
\begin{align}\label{42}
\Phi_4&=-i\int_{|\xi'|=1}\int^{+\infty}_{-\infty}{\rm trace} [\pi^+_{\xi_n}\sigma_{0}(\nabla_X^{S(TM)}\nabla_Y^{S(TM)}D^{-2})\times
\partial_{\xi_n}\sigma_{-(n-1)}(D^{-(n-2)})](x_0)d\xi_n\sigma(\xi')dx'\nonumber\\
&=i\int_{|\xi'|=1}\int^{+\infty}_{-\infty}{\rm trace} [\partial_{\xi_n}\pi^+_{\xi_n}\sigma_{0}(\nabla_X^{S(TM)}\nabla_Y^{S(TM)}D^{-2})\times
\sigma_{-(n-1)}(D^{-(n-2)})](x_0)d\xi_n\sigma(\xi')dx'.
\end{align}
By (3.30) in \cite{Wa6}, we have
\begin{align}\label{43}
\sigma_{-(n-1)}(D^{-n+2})
&=\frac{n-2}{2} (1+\xi_{n}^{2})^{(-\frac{n}{2}+2)}
   \Big[\frac{-i}{(1+\xi_n^2)^2}\times\frac{n+1}{2}h'(0)\xi_n-\frac{2ih'(0)\xi_n}{(1+\xi_n^2)^3} \Big]\nonumber\\
   &+\sqrt{-1} h'(0)(-\frac{n^{2}}{4}+\frac{3n}{2}-2)  \xi_{n} (1+\xi_{n}^{2})^{(-\frac{n}{2}-1)}.
\end{align}
 By Lemma \ref{lem3}, we have
\begin{align}\label{45}
\partial_{\xi_n}\pi^+_{\xi_n}\sigma_{0}(\nabla_X^{S(TM)}\nabla_Y^{S(TM)}D^{-2})&=-\frac{i}{2(\xi_n-i)^2}\sum_{j,l=1}^{n-1}X_jY_l\xi_j\xi_l-\frac{1}{2(\xi_n-i)^2}X_nY_n\nonumber\\
&+\frac{1}{2(\xi_n-i)^2}\sum_{j=1}^{n-1}(X_jY_n+X_nY_j)\xi_j.
\end{align}
We note that $i<n,~\int_{|\xi'|=1}\xi_{i_{1}}\xi_{i_{2}}\cdots\xi_{i_{2d+1}}\sigma(\xi')=0$,
so we omit some items that have no contribution for computing {\bf case b)}.\\
Then, we have
\begin{align}\label{39}
\Phi_4&=i\int_{|\xi'|}\int_{|\xi'|=1}\int^{+\infty}_{-\infty}{\rm tr}
\bigg\{\bigg(-\frac{i}{2(\xi_n-i)^2}\sum_{j,l=1}^{n-1}X_jY_l\xi_j\xi_l-\frac{1}{2(\xi_n-i)^2}X_nY_n\bigg)\times \bigg[\frac{n-2}{2} (1+\xi_{n}^{2})^{(-\frac{n}{2}+2)}
   \nonumber\\
   &\bigg(\frac{-i}{(1+\xi_n^2)^2}\times\frac{n+1}{2}h'(0)\xi_n-\frac{2ih'(0)\xi_n}{(1+\xi_n^2)^3} \bigg)+\sqrt{-1} h'(0)(-\frac{n^{2}}{4}
   +\frac{3n}{2}-2)  \xi_{n} (1+\xi_{n}^{2})^{(-\frac{n}{2}-1)}\bigg]
\bigg\}d\xi_n\sigma(\xi')dx' \nonumber\\
&=\frac{i h'(0)}{8} Vol(S^{n-2})\bigg(\frac{1}{n-1}\sum_{j=1}^{n-1}X_jY_j-iX_nY_n\bigg)\int_{\Gamma^+}\frac{-(n-2)(n+1)\xi_{n}^{3}+(-2n^{2}
   +3n-2)  \xi_{n}}{(\xi_n-i)^{(\frac{n}{2}+3)}(\xi_n+i)^{(\frac{n}{2}+1)}}d\xi_ndx'\nonumber\\
&= Vol(S^{n-2})\frac{2^{\frac{n}{2}}h'(0)\pi i}{4(\frac{n}{2}+2)!}\bigg(\frac{i}{n-1}g(X^T,Y^T)+X_nY_n\bigg)\left[\frac{-(n^2-n-2)\xi_{n}^{3}-(2n^{2}
-3n+2)\xi_{n}}{(\xi_n+i)^{\frac{n}{2}+1}}\right]^{(\frac{n}{2}+2)}\bigg|_{\xi_n=i}dx'\nonumber\\
   &:=Vol(S^{n-2})h'(0)2^{\frac{n}{2}}\bigg(\frac{i}{n-1}g(X^T,Y^T)+X_nY_n\bigg)\frac{\pi i}{4(\frac{n}{2}+2)!}C_{0}dx', 
\end{align}
where
\begin{align}
C_0&=\left[\frac{-(n-2)(n+1)\xi_{n}^{3}+(-2n^{2}
   +3n-2)  \xi_{n}}{(\xi_n+i)^{(\frac{n}{2}+1)}}\right]^{(\frac{n}{2}+2)}\bigg|_{\xi_n=i}\nonumber\\
 &=-i^{-n-4}2^{-n-1}\bigg((2n^2-2n-4)C_{-\frac{n}{2}-1}^{\frac{n}{2}-1}+(3n^2-3n-6)C_{-\frac{n}{2}-1}^{\frac{n}{2}}\nonumber\\
 &+(n^2-4)C_{-\frac{n}{2}-1}^{\frac{n}{2}+1}+(-\frac{n^2}{4}-\frac{n}{2}-1)C_{-\frac{n}{2}-1}^{\frac{n}{2}+2}\bigg)(\frac{n}{2}+2)!.\nonumber\\
\end{align}
\noindent {\bf  case c)}~$r=-1,~\ell=-(n-2),~k=j=|\alpha|=0$.\\
By (\ref{qqq15}), we get
\begin{align}\label{61}
\Phi_5=-i\int_{|\xi'|=1}\int^{+\infty}_{-\infty}{\rm trace} [\pi^+_{\xi_n}\sigma_{-1}(\nabla_X^{S(TM)}\nabla_Y^{S(TM)}D^{-2})\times
\partial_{\xi_n}\sigma_{-(n-2)}(D^{-(n-2)})](x_0)d\xi_n\sigma(\xi')dx'.
\end{align}
By (3.33) in \cite{Wa6}, we have
\begin{equation}\label{62}
\partial_{\xi_n}\sigma_{-(n-2)}(D^{-(n-2)})(x_0)
=\partial_{\xi_n}((|\xi|^2)^{1-\frac{n}{2}})(x_0)=2(1-\frac{n}{2})\xi_n(1+\xi_n^2)^{-\frac{n}{2}}.
\end{equation}
By Lemma \ref{lem356}, we have
\begin{align}\label{63}
\sigma_{-1}(\nabla_X^{S(TM)}\nabla_Y^{S(TM)}D^{-2})&=i\sum_jA(Y)X_j\xi_j|\xi|^{-2}+i\sum_lA(X)Y_l\xi_l|\xi|^{-2}+i\sum_{j,l=1}^{n-1}X_j\frac{\partial_{Y_l}}{\partial_{X_j}}\xi_l|\xi|^{-2}\nonumber\\
&+i\sum_{j,l=1}^{n-1}X_jY_l\xi_j\xi_l|\xi|^{-4}\xi_k(\Gamma^k-2\delta^k)+2i\sum_{j,l=1}^{n-1}X_jY_l\xi_j\xi_l|\xi|^{-6}\xi^j\xi_\alpha\xi_\beta\partial_jg^{\alpha\beta}\nonumber\\
&-\sum_{j=1}^{n-1}[\sum_{l=1}^{n-1}(X_jY_l+X_lY_j)\xi_l]-i\partial_{x_j}(|\xi|^{-2}).
\end{align}
By the Cauchy integral formula, we obtain
\begin{eqnarray}\label{64}
\pi^+_{\xi_n} \Big( \frac{\xi_n}{(1+\xi_n^2)^{2}}\Big)
&=&\frac{1}{2\pi i} \int_{\Gamma^+}\frac{\eta_n}{(\xi_n-\eta_n)(1+\eta_n^{2} )^2}d\eta_n\nonumber\\
&=&\left[\frac{\eta_n}{(\xi_n-\eta_n)(\eta_n^{2}+i)^{2}}\right]^{1}\Big|_{\eta_n=i}
     = \frac{-i}{4(\xi_n-i)^2},\\
\pi^+_{\xi_n} \Big( \frac{\xi_n}{(1+\xi_n^2)^{3}}\Big)
&=&\frac{1}{2\pi i} \int_{\Gamma^+}\frac{\eta_n}{(\xi_n-\eta_n)(1+\eta_n^{2} )^3}d\eta_n\nonumber\\
&=&\frac{1}{2} \left[\frac{\eta_n}{(\xi_n-\eta_n)(\eta_n^{2}+i )^{3}}\right]^{2}\Big|_{\eta_n=i}
     = \frac{-i}{16(\xi_n-i)^2}- \frac{1}{8(\xi_n-i)^3}.
\end{eqnarray}
We note that $i<n,~\int_{|\xi'|=1}\xi_{i_{1}}\xi_{i_{2}}\cdots\xi_{i_{2d+1}}\sigma(\xi')=0$,
so we omit some items that have no contribution for computing {\bf case c)}, then
 \begin{align}\label{65}
\pi^+_{\xi_n}\sigma_{-1}(\nabla_X^{S(TM)}\nabla_Y^{S(TM)}D^{-2})
         &=\sum_{j=1}^{n-1}X_jY_j\xi_j\xi_lh'(0)\Big(\frac{1}{8(\xi_n-i)^2}-\frac{i}{4(\xi_n-i)^3} \Big)+2iX_nY_nh'(0)\xi_n.
\end{align}
Therefore
\begin{align}\label{66}
&\int_{|\xi'|=1}\int^{+\infty}_{-\infty}{\rm tr}
\bigg[\sum_{j=1}^{n-1}X_jY_j\xi_j\xi_lh'(0)\bigg(\frac{1}{8(\xi_n-i)^2}-\frac{i}{4(\xi_n-i)^3} \bigg)
\times\frac{2(1-\frac{n}{2})\xi_n}{(1+\xi_n^2)^{\frac{n}{2}}}
\bigg]d\xi_n\sigma(\xi')dx'\nonumber\\
&=\int_{|\xi'|=1}\frac{1}{4}\sum_{j=1}^{n-1}X_jY_j\xi_j\xi_lh'(0)2(1-\frac{n}{2})2^{\frac{n}{2}}\sigma(\xi')\int_{\Gamma^+}\frac{\xi_n}{(\xi_n-i)^{\frac{n}{2}+2}(\xi_n+i)^{\frac{n}{2}}}dx'\nonumber\\
&+\int_{|\xi'|=1}\frac{1}{2}\sum_{j=1}^{n-1}X_jY_j\xi_j\xi_lh'(0)2(1-\frac{n}{2})2^{\frac{n}{2}}\sigma(\xi')\int_{\Gamma^+}\frac{\xi_n}{(\xi_n-i)^{\frac{n}{2}+3}(\xi_n+i)^{\frac{n}{2}}}dx'\nonumber\\
&=Vol(S^{n-2})\frac{1}{4}g(X^T,Y^T)h'(0)2(1-\frac{n}{2})2^{\frac{n}{2}}\frac{2\pi i}{(1+\frac{n}{2})!}\left[\frac{\xi_{n}}{(\xi_n+i)^{\frac{n}{2}}}\right]^{(\frac{n}{2}+1)}\bigg|_{\xi_n=i}dx'\nonumber\\
&-Vol(S^{n-2})\frac{1}{2}g(X^T,Y^T)h'(0)2(1-\frac{n}{2})2^{\frac{n}{2}}\frac{2\pi i}{(2+\frac{n}{2})!}\left[\frac{\xi_{n}}{(\xi_n+i)^{\frac{n}{2}}}\right]^{(\frac{n}{2}+2)}\bigg|_{\xi_n=i}dx'\nonumber\\
&:=Vol(S^{n-2})g(X^T,Y^T)h'(0)\frac{1}{n-1}(1-\frac{n}{2})2^{\frac{n}{2}}\frac{\pi i}{(1+\frac{n}{2})!}D_0dx'\nonumber\\
&-Vol(S^{n-2})g(X^T,Y^T)h'(0)\frac{1}{n-1}(1-\frac{n}{2})2^{\frac{n}{2}}\frac{2\pi i}{(2+\frac{n}{2})!}D_1dx'\nonumber\\
\end{align}
\begin{align}\label{6IIII}
&\int_{|\xi'|=1}\int^{+\infty}_{-\infty}{\rm tr}
\bigg[2iX_nY_nh'(0)\xi_n
\times\frac{2(1-\frac{n}{2})\xi_n}{(1+\xi_n^2)^{\frac{n}{2}}}
\bigg]d\xi_n\sigma(\xi')dx'\nonumber\\
&=\int_{|\xi'|=1}4ih'(0)X_nY_n(1-\frac{n}{2})2^{\frac{n}{2}}\int_{\Gamma^+}\frac{\xi_n^2}{(\xi_n-i)^{\frac{n}{2}+2}(\xi_n+i)^{\frac{n}{2}}}dx'\nonumber\\
&=-Vol(S^{n-2})h'(0)X_nY_n(1-\frac{n}{2})2^{\frac{n}{2}}\frac{16\pi}{(\frac{n}{2}+1)!}\left[\frac{\xi_{n}^2}{(\xi_n+i)^{\frac{n}{2}}}\right]^{(\frac{n}{2}+1)}\bigg|_{\xi_n=i}dx'\nonumber\\
&:=-Vol(S^{n-2})h'(0)X_nY_n(1-\frac{n}{2})2^{\frac{n}{2}}\frac{16\pi}{(\frac{n}{2}+1)!}D_2dx'\nonumber\\
\end{align}
Then
\begin{align}\label{6666}
\Phi_5&=-iVol(S^{n-2})h'(0)(1-\frac{n}{2})2^{\frac{n}{2}}\bigg\{\frac{1}{n-1}g(X^T,Y^T)\bigg(\frac{\pi i}{(1+\frac{n}{2})!}D_0-\frac{2\pi i}{(2+\frac{n}{2})!}D_1\bigg)-X_nY_n\frac{16\pi}{(\frac{n}{2}+1)!}D_2\bigg\}dx',
\end{align}
where
\begin{align}
D_0&=\left[\frac{\xi_{n}}{(\xi_n+i)^{\frac{n}{2}}}\right]^{(\frac{n}{2}+1)}\bigg|_{\xi_n=i}=-i^{-n-2}2^{-n-1}\bigg(2C_{-\frac{n}{2}}^{\frac{n}{2}}+C_{-\frac{n}{2}}^{\frac{n}{2}+1}\bigg)(\frac{n}{2}+1)!;\nonumber\\
D_1&=\left[\frac{\xi_{n}}{(\xi_n+i)^{\frac{n}{2}}}\right]^{(\frac{n}{2}+2)}\bigg|_{\xi_n=i}=(2i)^{-n-2}\bigg(2iC_{-\frac{n}{2}}^{\frac{n}{2}+1}+iC_{-\frac{n}{2}}^{\frac{n}{2}+2}\bigg)(\frac{n}{2}+2)!;\nonumber\\
D_2&=\left[\frac{\xi_{n}^2}{(\xi_n+i)^{\frac{n}{2}}}\right]^{(\frac{n}{2}+1)}\bigg|_{\xi_n=i}=-(2i)^{-n-1}\bigg(4C_{-\frac{n}{2}}^{\frac{n}{2}-1}+4C_{-\frac{n}{2}}^{\frac{n}{2}}+C_{-\frac{n}{2}}^{\frac{n}{2}+1}\bigg)(\frac{n}{2}+1)!.\nonumber\\
\end{align}
Now $\Phi$ is the sum of the cases (a), (b) and (c). Therefore, we get
\begin{align}\label{795}
\Phi&=\sum_{i=1}^5\Phi_i\nonumber\\
&=-\frac{1}{2}Vol(S^{n-2})h'(0)2^{\frac{n}{2}}\bigg\{(\frac{n}{2}-1)\bigg(\frac{1}{n-1}g(X^T,Y^T)+X_nY_n\bigg)\frac{\pi}{2(\frac{n}{2}+2)!}A_{0}+(\frac{n}{2}-1)X_nY_n\frac{\pi i}{(\frac{n}{2}+2)!}A_{1}\nonumber\\
&+(\frac{n}{2}-1)\bigg[\frac{i}{n-1}\bigg(h'(0)g(X^T,Y^T)+\partial_{x_n}(g(X^T,Y^T))\bigg)-\partial_{x_n}(X_nY_n)\bigg]\frac{\pi i}{(\frac{n}{2}+1)!}A_2-(1-\frac{n}{2})\nonumber\\
&\bigg(\frac{1}{n-1}g(X^T,Y^T)-X_nY_n\bigg)\frac{2\pi i}{(\frac{n}{2}+2)!}B_0-\bigg(\frac{i}{n-1}g(X^T,Y^T)+X_nY_n\bigg)\frac{\pi i}{2(\frac{n}{2}+2)!}C_{0}-\frac{1}{n-1}g(X^T,Y^T)\nonumber\\
&(1-\frac{n}{2})\frac{2\pi }{(1+\frac{n}{2})!}D_0+\frac{1}{n-1}g(X^T,Y^T)(1-\frac{n}{2})\frac{4\pi }{(2+\frac{n}{2})!}D_1-X_nY_n(1-\frac{n}{2})\frac{32\pi i }{(\frac{n}{2}+1)!}D_2\bigg\}dx'.\nonumber\\
\end{align}
Then, by (\ref{795}), we obtain following theorem
\begin{thm}\label{thmb1}
Let $M$ be an $n$-dimensional oriented
compact spin manifold with boundary $\partial M$ and n is even, then we get the following equality:
\begin{align}
\label{b2773}
&\widetilde{{\rm Wres}}[\pi^+(\nabla_X^{S(TM)}\nabla_Y^{S(TM)}D^{-2})\circ\pi^+(D^{-(n-2)})]\nonumber\\
&=\frac{(2\pi)^{\frac{n}{2}}}{3(\frac{n}{2}-1)!}\int_{M}[Ric(X,Y)-\frac{1}{2}sg(X,Y)]d{\rm Vol_{M}}+\frac{(2\pi)^{\frac{n}{2}}}{4(\frac{n}{2}-1)!}\int_{M}sg(X,Y)d{\rm Vol_{M}}-\int_{\partial M}\frac{1}{2}Vol(S^{n-2})h'(0)2^{\frac{n}{2}}\nonumber\\
&\bigg\{(\frac{n}{2}-1)\bigg(\frac{1}{n-1}g(X^T,Y^T)+X_nY_n\bigg)\frac{\pi}{2(\frac{n}{2}+2)!}A_{0}+(\frac{n}{2}-1)X_nY_n\frac{\pi i}{(\frac{n}{2}+2)!}A_{1}+(\frac{n}{2}-1)\bigg[\frac{i}{n-1}\bigg(h'(0)\nonumber\\
&g(X^T,Y^T)+\partial_{x_n}(g(X^T,Y^T))\bigg)-\partial_{x_n}(X_nY_n)\bigg]\frac{\pi i}{(\frac{n}{2}+1)!}A_2-(1-\frac{n}{2})\bigg(\frac{1}{n-1}g(X^T,Y^T)-X_nY_n\bigg)\nonumber\\
&\frac{2\pi i}{(\frac{n}{2}+2)!}B_0-\bigg(\frac{i}{n-1}g(X^T,Y^T)+X_nY_n\bigg)\frac{\pi i}{2(\frac{n}{2}+2)!}C_{0}-\frac{1}{n-1}g(X^T,Y^T)(1-\frac{n}{2})\frac{2\pi }{(1+\frac{n}{2})!}D_0\nonumber\\
&+\frac{1}{n-1}g(X^T,Y^T)(1-\frac{n}{2})\frac{4\pi }{(2+\frac{n}{2})!}D_1-X_nY_n(1-\frac{n}{2})\frac{32\pi i }{(\frac{n}{2}+1)!}D_2\bigg\}d{\rm Vol_{M}}.\nonumber\\
\end{align}
\end{thm}
Next, we compute the noncommutative residue $\widetilde{{\rm Wres}}[\pi^+(\nabla_X^{S(TM)}\nabla_Y^{S(TM)}D^{-1})\circ\pi^+(D^{-(n-1)})]$ on $n$-dimensional oriented compact manifolds with boundary, and $n$ is even.\\
\indent Similar to \cite{Wa3},  then we can compute the noncommutative residue
\begin{align}
\label{c1}
&\widetilde{{\rm Wres}}[\pi^+(\nabla_X^{S(TM)}\nabla_Y^{S(TM)}D^{-1})\circ\pi^+(D^{-(n-1)})]\nonumber\\
&=\int_M\int_{|\xi|=1}{\rm
trace}_{\wedge^*T^*M\bigotimes\mathbb{C}}[\sigma_{-n}(\nabla_X^{S(TM)}\nabla_Y^{S(TM)}D^{-1}\circ D^{-(n-1)})]\sigma(\xi)dx+\int_{\partial M}\widetilde{\Phi},
\end{align}
where
\begin{align}
\label{c2}
 \widetilde{\Phi}&=\int_{|\xi'|=1}\int^{+\infty}_{-\infty}\sum^{\infty}_{j, k=0}\sum\frac{(-i)^{|\alpha|+j+k+1}}{\alpha!(j+k+1)!}
\times {\rm trace}_{\wedge^*T^*M\bigotimes\mathbb{C}}[\partial^j_{x_n}\partial^\alpha_{\xi'}\partial^k_{\xi_n}\sigma^+_{r}(\nabla_X^{S(TM)}\nabla_Y^{S(TM)}D^{-1})(x',0,\xi',\xi_n)
\nonumber\\
&\times\partial^\alpha_{x'}\partial^{j+1}_{\xi_n}\partial^k_{x_n}\sigma_{l}(D^{-(n-1)})(x',0,\xi',\xi_n)]d\xi_n\sigma(\xi')dx',
\end{align}
and the sum is taken over $r+l-k-j-|\alpha|-1=-n,~~r\leq 1,~~l\leq -(n-1)$.\\

\indent Similarly, the interior of $\widetilde{{\rm Wres}}[\pi^+(\nabla_X^{S(TM)}\nabla_Y^{S(TM)}D^{-1})\circ\pi^+(D^{-(n-1)})]$ has been given in Theorem \ref{thm2}. Therefore, we  need to compute $\int_{\partial M} \widetilde{\Phi}$. When sum is taken over $
r+l-k-j-|\alpha|=-n+1,~~r\leq 1,~~l\leq -(n-1),$  we have the following five cases:
~\\
\noindent  {\bf case a)~I)}~$r=1,~l=-(n-1),~k=j=0,~|\alpha|=1$.\\
\noindent By (\ref{c2}), we get
\begin{equation}
\label{c4}
\widetilde{\Phi}_1=-\int_{|\xi'|=1}\int^{+\infty}_{-\infty}\sum_{|\alpha|=1}
{\rm trace}[\partial^\alpha_{\xi'}\pi^+_{\xi_n}\sigma_{1}(\nabla_X^{S(TM)}\nabla_Y^{S(TM)}D^{-1})\times
 \partial^\alpha_{x'}\partial_{\xi_n}\sigma_{-(n-1)}(D^{-(n-1)})](x_0)d\xi_n\sigma(\xi')dx'.
\end{equation}
By Lemma 2.2 in \cite{Wa3}, for $i<n$, then
\begin{align}
\label{c5}
\partial_{x_i}\sigma_{-(n-1)}(D^{-(n-1)})(x_0)&=\partial_{x_i}{(ic(\xi)|\xi|^{-n})}(x_0)\nonumber\\
&=i\partial_{x_i}c(\xi)(x_0)|\xi|^{-n}+ic(\xi)\partial_{x_i}(|\xi|^{-n})(x_0)
 \partial_{x_i}(|\xi|^2)(x_0)\nonumber\\
 &=0,
\end{align}
\noindent so $\widetilde{\Phi}_1=0$.\\
 \noindent  {\bf case a)~II)}~$r=1,~l=-(n-1),~k=|\alpha|=0,~j=1$.\\
\noindent By (\ref{c2}), we get
\begin{align}
\label{c6}
\widetilde{\Phi}_2&=-\frac{1}{2}\int_{|\xi'|=1}\int^{+\infty}_{-\infty} {\rm
trace} [\partial_{x_n}\pi^+_{\xi_n}\sigma_{1}(\nabla_X^{S(TM)}\nabla_Y^{S(TM)}D^{-1})\times
\partial_{\xi_n}^2\sigma_{-(n-1)}(D^{-(n-1)})](x_0)d\xi_n\sigma(\xi')dx'\nonumber\\
&=-\frac{1}{2}\int_{|\xi'|=1}\int^{+\infty}_{-\infty} {\rm
trace} [\partial_{\xi_n}^2\partial_{x_n}\pi^+_{\xi_n}\sigma_{1}(\nabla_X^{S(TM)}\nabla_Y^{S(TM)}D^{-1})\times
\sigma_{-(n-1)}(D^{-(n-1)})](x_0)d\xi_n\sigma(\xi')dx'.
\end{align}
\noindent By (3.11) in \cite{Wa7}, we have\\
\begin{eqnarray}\label{c7}
\sigma_{-(n-1)}(D^{-(n-1)})
&=&\frac{\sqrt{-1}[c(\xi')+\xi_nc(\mathrm{d}x_n)]}{(1+\xi_n^2)^{\frac{n}{2}}}.
\end{eqnarray}
By (\ref{mki}), we have
\begin{align}\label{c8}
&\partial_{x_n}\sigma_{1}(\nabla_X^{S(TM)}\nabla_Y^{S(TM)}D^{-1})\nonumber\\
&=\partial_{x_n}\bigg(-i\sum_{j,l=1}^nX_jY_l\xi_j\xi_l\frac{c(\xi)}{|\xi|^{2}}\bigg)\nonumber\\
&=-i\sum_{j,l=1}^n\left(\frac{\partial{X_j}}{\partial_{x_n}}Y_l\xi_j\xi_l+X_j\frac{\partial{Y_l}}{\partial_{x_n}}\xi_j\xi_l\right)\frac{c(\xi)}{|\xi|^{2}}+\sum_{j,l=1}^nX_jY_l\xi_j\xi_l \left(\frac{\partial_{x_n}c(\xi')}{1+\xi_n^2}+\frac{c(\xi)h'(0)|\xi'|^2}{(1+\xi_n^2)^2}\right).
\end{align}
We note that $i<n,~\int_{|\xi'|=1}\xi_{i_{1}}\xi_{i_{2}}\cdots\xi_{i_{2d+1}}\sigma(\xi')=0$,
so we omit some items that have no contribution for computing {\bf case a)~II)}. Then, we have
\begin{align}\label{c9}
&\partial_{x_n}\pi^+_{\xi_n}\sigma_{1}(\nabla_X^{S(TM)}\nabla_Y^{S(TM)}D^{-1})\nonumber\\
&=\pi^+_{\xi_n}\partial_{x_n}\sigma_{1}(\nabla_X^{S(TM)}\nabla_Y^{S(TM)}D^{-1})\nonumber\\
&=\bigg(\sum_{j,l=1}^{n-1}\frac{\partial{X_j}}{\partial_{x_n}}Y_l\xi_j\xi_l+\frac{\partial{X_n}}{\partial_{x_n}}Y_n\xi_n^2+\sum_{j,l=1}^{n-1}X_j\frac{\partial{Y_l}}{\partial_{x_n}}\xi_j\xi_l+X_n\frac{\partial{Y_n}}{\partial_{x_n}}\xi_n^2\bigg)\times-i\frac{c(\xi')+ic(\mathrm{d}x_n)}{2(\xi_n-i)}\nonumber\\
&+i\sum_{j,l=1}^{n-1}X_jY_l\xi_j\xi_lh'(0)|\xi'|^2\left(\frac{ic(\xi')}{4(\xi_n-i)}+\frac{c(\xi')+ic(\mathrm{d}x_n)}{4(\xi_n-i)^2}\right)-\sum_{j,l=1}^{n-1}X_jY_l\xi_j\xi_l\frac{\partial_{x_n}c(\xi')}{2(\xi_n-i)}\nonumber\\
&-iX_nY_n\left\{\frac{\partial_{x_n}c(\xi')}{2(\xi_n-i)}
 +h'(0)|\xi'|\left[-\frac{2ic(\xi')-3c(\mathrm{d}x_n)}{4(\xi_n-i)}+\frac{[c(\xi')+ic(\mathrm{d}x_n)][i(\xi_n-i)+1]}{4(\xi_n-i)^2}\right]\right\}.\nonumber\\
\end{align}
By further calculation, we have
\begin{align}\label{555}
&\partial_{\xi_n}^2\partial_{x_n}\pi^+_{\xi_n}\sigma_{1}(\nabla_X^{S(TM)}\nabla_Y^{S(TM)}D^{-1})\nonumber\\
&=\sum_{j,l=1}^{n-1}\left(\frac{\partial{X_j}}{\partial_{x_n}}Y_l+X_j\frac{\partial{Y_l}}{\partial_{x_n}}\right)\xi_j\xi_l\times\frac{-i[c(\xi')+ic(\mathrm{d}x_n)]}{(\xi_n-i)^3}\left(\frac{\partial{X_n}}{\partial_{x_n}}Y_n+X_n\frac{\partial{Y_n}}{\partial_{x_n}}\right)\times\frac{-i(2+\xi_ni)[c(\xi')+ic(\mathrm{d}x_n)]}{(\xi_n-i)^3}\nonumber\\
&+i\sum_{j,l=1}^{n-1}X_jY_l\xi_j\xi_lh'(0)|\xi'|^2\left[\frac{ic(\xi')}{2(\xi_n-i)^3}+\frac{3[c(\xi')+ic(\mathrm{d}x_n)]}{2(\xi_n-i)^4}\right]-\sum_{j,l=1}^{n-1}X_jY_l\xi_j\xi_l\frac{\partial_{x_n}c(\xi')}{(\xi_n-i)^3}\nonumber\\
&-iX_nY_n\left\{\frac{\partial_{x_n}c(\xi')}{(\xi_n-i)^3}
 +h'(0)|\xi'|\left[-\frac{2ic(\xi')-3c(\mathrm{d}x_n)}{2(\xi_n-i)^3}+\frac{[ic(\xi')-c(\mathrm{d}x_n)+ic(\mathrm{d}x_n)+c(\xi')](\xi_n+2i)}{2(\xi_n-i)^4}\right]\right\}.\nonumber\\
\end{align}
Moreover, we have
\begin{align}\label{c10}
&\int_{|\xi'|}\int_{-\infty}^{\infty}{\rm tr}\bigg[\sum_{j,l=1}^{n-1}\left(\frac{\partial{X_j}}{\partial_{x_n}}Y_l+X_j\frac{\partial{Y_l}}{\partial_{x_n}}\right)\xi_j\xi_l\frac{-i[c(\xi')+ic(\mathrm{d}x_n)]}{(\xi_n-i)^3}\times
 \frac{\sqrt{-1}[c(\xi')+\xi_nc(\mathrm{d}x_n)]}{(1+\xi_n^2)^{\frac{n}{2}}}\bigg]d\xi_n\sigma(\xi')dx'\nonumber\\
 &=\int_{|\xi'|}\int_{-\infty}^{\infty}{\rm tr}\bigg[\sum_{j,l=1}^{n-1}\left(\frac{\partial{X_j}}{\partial_{x_n}}Y_l+X_j\frac{\partial{Y_l}}{\partial_{x_n}}\right)\xi_j\xi_l\frac{i\xi_n-1}{(\xi_n+i)^{\frac{n}{2}}(\xi_n-i)^{\frac{n}{2}+3}}\bigg]d\xi_n\sigma(\xi')dx'\nonumber\\
 &=\int_{|\xi'|}\sum_{j,l=1}^{n-1}\left(\frac{\partial{X_j}}{\partial_{x_n}}Y_l+X_j\frac{\partial{Y_l}}{\partial_{x_n}}\right)\xi_j\xi_l\sigma(\xi')2^{\frac{n}{2}} \frac{2\pi i}{(\frac{n}{2}+2)!}
\left[\frac{i\xi_n-1}{(\xi_n+i)^{\frac{n}{2}}}\right]^{(\frac{n}{2}+2)}
\bigg|_{\xi_n=i}dx'\nonumber\\
 &=Vol(S^{n-2})\frac{1}{n-1}\bigg(h'(0)g(X^T,Y^T)+\partial_{x_n}(g(X^T,Y^T))\bigg)2^{\frac{n}{2}} \frac{2\pi i}{(\frac{n}{2}+2)!}
\left[\frac{i\xi_n-1}{(\xi_n+i)^{\frac{n}{2}}}\right]^{(\frac{n}{2}+2)}
\bigg|_{\xi_n=i}dx'\nonumber\\
 &:=Vol(S^{n-2})\frac{1}{n-1}\bigg(h'(0)g(X^T,Y^T)+\partial_{x_n}(g(X^T,Y^T))\bigg)2^{\frac{n}{2}} \frac{2\pi i}{(\frac{n}{2}+2)!}E_{0}dx'.
\end{align}
Similarly, we can obtain items of other items. Therefore, we get\\
\begin{align}\label{c11}
\widetilde{\Phi}_2
&=-\frac{1}{2}Vol(S^{n-2})2^{\frac{n}{2}}\bigg\{\frac{1}{n-1}\bigg(h'(0)g(X^T,Y^T)+\partial_{x_n}(g(X^T,Y^T))\bigg)\frac{2\pi i}{(\frac{n}{2}+2)!}E_{0}+\partial_{x_n}(X_nY_n)\nonumber\\
&\frac{2\pi i}{(\frac{n}{2}+2)!}E_{1}+(i-1)h'(0)\bigg(\frac{1}{n-1}g(X^T,Y^T)-X_nY_n\bigg)\frac{\pi i}{(\frac{n}{2}+2)!}E_{2}-\frac{1}{n-1}g(X^T,Y^T)h'(0)\frac{\pi i}{(\frac{n}{2}+3)!}E_{3}\nonumber\\
&+X_nY_n\frac{\pi i}{(\frac{n}{2}+3)!}E_{4}\bigg\}dx',
\end{align}
where
\begin{align}
E_0&=\left[\frac{i\xi_n-1}{(\xi_n+i)^{\frac{n}{2}}}\right]^{(\frac{n}{2}+2)}
\bigg|_{\xi_n=i}=-i^{-n-4}2^{-n-1}\bigg(C_{-\frac{n}{2}}^{\frac{n}{2}+1}+C_{-\frac{n}{2}}^{\frac{n}{2}+2}\bigg)(\frac{n}{2}+2)!;\nonumber\\
E_1&=\left[\frac{-\xi_n^2+i\xi_n-2}{(\xi_n+i)^{\frac{n}{2}}}\right]^{(\frac{n}{2}+2)}
\bigg|_{\xi_n=i}=-i^{-n-4}2^{-n-1}\bigg(2C_{-\frac{n}{2}}^{\frac{n}{2}}+C_{-\frac{n}{2}}^{\frac{n}{2}+1}-C_{-\frac{n}{2}}^{\frac{n}{2}+2}\bigg)(\frac{n}{2}+2)!;\nonumber\\
E_2&=\left[\frac{1}{(\xi_n+i)^{\frac{n}{2}}}\right]^{(\frac{n}{2}+2)}
\bigg|_{\xi_n=i}=(2i)^{-n-2}A_{-\frac{n}{2}}^{\frac{n}{2}+2}(\frac{n}{2}+2)!=B_0;\nonumber\\
E_3&=\left[\frac{i\xi_n-1}{(\xi_n+i)^{\frac{n}{2}}}\right]^{(\frac{n}{2}+3)}
\bigg|_{\xi_n=i}=-i^{-n-5}2^{-n-2}\bigg(C_{-\frac{n}{2}}^{\frac{n}{2}+2}+C_{-\frac{n}{2}}^{\frac{n}{2}+3}\bigg)(\frac{n}{2}+3)!;\nonumber\\
E_4&=\left[\frac{-(2+i)\xi_n^2+(4+i)\xi_n-2i}{(\xi_n+i)^{\frac{n}{2}}}\right]^{(\frac{n}{2}+3)}
\bigg|_{\xi_n=i}\nonumber\\
&=-i^{-n-6}2^{-\frac{n}{2}-1}\bigg((8+4i)C_{-\frac{n}{2}}^{\frac{n}{2}}+(12+6i)C_{-\frac{n}{2}}^{\frac{n}{2}+1}+(10+4i)C_{-\frac{n}{2}}^{\frac{n}{2}+2}+(2+i)C_{-\frac{n}{2}}^{\frac{n}{2}+3}\bigg)(\frac{n}{2}+3)!.\nonumber\\
\end{align}
\noindent  {\bf case a)~III)}~$r=1,~l=-(n-1),~j=|\alpha|=0,~k=1$.\\
\noindent By (\ref{c2}), we get
\begin{align}\label{c12}
\widetilde{\Phi}_3&=-\frac{1}{2}\int_{|\xi'|=1}\int^{+\infty}_{-\infty}
{\rm trace} [\partial_{\xi_n}\pi^+_{\xi_n}\sigma_{1}(\nabla_X^{S(TM)}\nabla_Y^{S(TM)}D^{-1})\times
\partial_{\xi_n}\partial_{x_n}\sigma_{-(n-1)}(D^{-(n-1)})](x_0)d\xi_n\sigma(\xi')dx'\nonumber\\
&=\frac{1}{2}\int_{|\xi'|=1}\int^{+\infty}_{-\infty}
{\rm trace} [\partial_{\xi_n}^2\pi^+_{\xi_n}\sigma_{1}(\nabla_X^{S(TM)}\nabla_Y^{S(TM)}D^{-1})\times
\partial_{x_n}\sigma_{-(n-1)}(D^{-(n-1)})](x_0)d\xi_n\sigma(\xi')dx'.
\end{align}
\noindent By (3.17) in \cite{Wa7}, we have
\begin{eqnarray}\label{c13}
\partial_{x_n} \big(\sigma_{-(n-1)}(D^{-(n-1)})\big)(x_0)
=\frac{i\partial_{x_n}c(\xi')(x_0)}{(1+\xi_n^2)^{\frac{n}{2}}}-\frac{\sqrt{-1}nh'(0)c(\xi)}{2(1+\xi_n^2)^{\frac{n}{2}+1}}.
\end{eqnarray}
By lemma \ref{lem9356}, we have 
\begin{align}\label{c14}
\partial_{\xi_n}^2\pi^+_{\xi_n}\sigma_{1}(\nabla_X^{S(TM)}\nabla_Y^{S(TM)}D^{-1})&=-\frac{c(\xi')+ic(\mathrm{d}x_n)}{(\xi_n-i)^3}\bigg(\sum_{j,l=1}^{n-1}X_jY_l\xi_j\xi_l
+X_nY_n\bigg).
\end{align}
Then, we have
\begin{align}\label{c16}
& \int_{|\xi'|}\int_{-\infty}^{\infty}{\rm tr}\bigg[-\frac{c(\xi')+ic(\mathrm{d}x_n)}{(\xi_n-i)^3}[\sum_{j,l=1}^{n-1}X_jY_l\xi_j\xi_l
+X_nY_n]\times\frac{i\partial_{x_n}c(\xi')(x_0)}{(1+\xi_n^2)^{\frac{n}{2}}}-\frac{\sqrt{-1}nh'(0)c(\xi)}{2(1+\xi_n^2)^{\frac{n}{2}+1}}\bigg](x_0)d\xi_n\sigma(\xi')dx'\nonumber\\
&= \int_{|\xi'|}\int_{-\infty}^{\infty}-2^{\frac{n}{2}}h'(0)\bigg(\sum_{j,l=1}^{n-1}X_jY_l\xi_j\xi_l
+X_nY_n\bigg)\frac{-i\xi_n^2-n\xi_n+ni-i}{(\xi_n+i)^{\frac{n}{2}+1}(\xi_n-i)^{\frac{n}{2}+4}}d\xi_n\sigma(\xi')dx'\nonumber\\\nonumber\\
 &=-Vol(S^{n-2})h'(0)2^{\frac{n}{2}+1}\bigg(\frac{1}{n-1}g(X^T,Y^T)+X_nY_n\bigg)\frac{\pi i}{(\frac{n}{2}+3)!}\left[\frac{-i\xi_n^2-n\xi_n+ni-i}{(\xi_n+i)^{(\frac{n}{2}+1)}}\right]
 ^{(\frac{n}{2}+3)}\bigg|_{\xi_n=i}dx'\nonumber\\
  &:=-Vol(S^{n-2})h'(0)2^{\frac{n}{2}+1}\bigg(\frac{1}{n-1}g(X^T,Y^T)+X_nY_n\bigg)\frac{\pi i}{(\frac{n}{2}+3)!}F_0dx'.
\end{align}
Therefore, we get
\begin{align}\label{c18}
\widetilde{\Phi}_3&=-Vol(S^{n-2})h'(0)2^{\frac{n}{2}+1}\bigg(\frac{1}{n-1}g(X^T,Y^T)+X_nY_n\bigg)\frac{\pi i}{(\frac{n}{2}+3)!}F_0dx',
\end{align}
where
\begin{align}
F_0&=\left[\frac{-i\xi_n^2-n\xi_n+ni-i}{(\xi_n+i)^{(\frac{n}{2}+1)}}\right]
 ^{(\frac{n}{2}+3)}\bigg|_{\xi_n=i}=i^{-n-6}2^{-n-1}\bigg(C_{-\frac{n}{2}}^{\frac{n}{2}+1}-(\frac{n}{2}-1)C_{-\frac{n}{2}}^{\frac{n}{2}+3}\bigg)(\frac{n}{2}+2)!.\nonumber\\
\end{align}
\noindent  {\bf case b)}~$r=0,~l=-(n-1),~k=j=|\alpha|=0$.\\
\noindent By (\ref{c2}), we get
\begin{align}\label{c19}
\widetilde{\Phi}_4&=-i\int_{|\xi'|=1}\int^{+\infty}_{-\infty}{\rm trace} [\pi^+_{\xi_n}\sigma_{0}(\nabla_X^{S(TM)}\nabla_Y^{S(TM)}D^{-1})\times
\partial_{\xi_n}\sigma_{-(n-1)}(D^{-(n-1)})](x_0)d\xi_n\sigma(\xi')dx'.
\end{align}
 By (3.24) in \cite{Wa7}, we have
\begin{align}\label{c20}
\partial_{\xi_n}\sigma_{-(n-1)}(D^{-(n-1)})(x_0)\bigg|_{|\xi'|=1}
&=\sqrt{-1}\bigg(\frac{c(\mathrm{d}x_n)}{(1+\xi_n^2)^{\frac{n}{2}}}-\frac{n[\xi_nc(\xi')+\xi_n^2c(\mathrm{d}x_n)]}{(1+\xi_n^2)^{\frac{n}{2}+1}} \bigg).
\end{align}
By Lemma \ref{lem356}, we have
\begin{align}\label{ggg}
\sigma_{0}(\nabla_X^{S(TM)}\nabla_Y^{S(TM)}D^{-1})&=-\sum_{j,l=1}^{n}X_jY_l\xi_j\xi_l\left[\frac{c(\xi)\sigma_0(D)c(\xi)}{|\xi|^4}
+\frac{c(\xi)}{|\xi|^6}\Sigma_jc(\mathrm{d}x_j)(\partial_{x_j}[c(\xi)]|\xi|^2-c(\xi)\partial_{x_j}(|\xi|^2))\right]\nonumber\\
&-i\sum_{j,l=1}^nX_j\frac{\partial_{Y_l}}{\partial_{x_j}}\xi_j\frac{c(\xi)}{|\xi|^2}
-\sum_jA(Y)X_j\xi_j\frac{c(\xi)}{|\xi|^2}-\sum_lA(Y)Y_l\xi_l\frac{c(\xi)}{|\xi|^2}\nonumber\\
&-\sum_{j}^{n}\sum_{l}^{n}(X_jY_l+X_lY_j)\xi_l)\partial_{x_j}\bigg(\frac{c(\xi)}{|\xi|^2}\bigg).\nonumber\\
\end{align}
\begin{align}\label{c21}
&\pi^+_{\xi_n}\sigma_{0}(\nabla_X^{S(TM)}\nabla_Y^{S(TM)}D^{-1})\nonumber\\
&=-\bigg(\sum_{j=1}^nX_j\frac{\partial{Y_n}}{\partial_{x_j}}+A(Y)X_n+A(X)Y_n\bigg)\frac{(\xi_n+i)c(\xi')-2c(\mathrm{d}x_n)}{4(\xi_n-i)^2}-\sum_{j=1}^{n-1}X_jY_l\xi_j\xi_l(Q_1-Q_2)\nonumber\\
&-X_nY_n(P_1+P_2)+\sum_{l=1}^{n-1}(X_nY_l+X_lY_n)\xi_l\bigg(\frac{i\partial_{x_n}(c(\xi'))}{2(\xi_n-i)}+h'(0)\frac{(-2-i\xi_n)c(\xi')}{4(\xi_n-i)^2}-h'(0)\frac{ic(\mathrm{d}x_n)}{4(\xi_n-i)^2}\bigg)\nonumber\\
&+X_nY_n\bigg(\frac{-i\partial_{x_n}(c(\xi'))}{(\xi_n-i)}+h'(0)\frac{(-i\xi_n)c(\xi')}{2(\xi_n-i)^2}-h'(0)\frac{-i\xi_nc(\mathrm{d}x_n)}{4(\xi_n-i)^2}\bigg),\nonumber\\
\end{align}
where
\begin{align}\label{m1}
q_1&=ic(\xi')p_0c(\xi')+ic(\mathrm{d}x_n)p_0c(\mathrm{d}x_n)+ic(\xi')c(\mathrm{d}x_n)\partial_{x_n}[c(\xi')];\nonumber\\
q_2&=[c(\xi')+ic(\mathrm{d}x_n)]p_0[c(\xi')+ic(\mathrm{d}x_n)]+c(\xi')c(\mathrm{d}x_n)\partial_{x_n}[c(\xi')]-i\partial_{x_n}[c(\xi')];\nonumber\\
Q_1&=-\frac{q_1}{4(\xi_n-i)}-\frac{q_2}{4(\xi_n-i)^2};\nonumber\\
Q_2&=\frac{h'(0)}{2}\left[\frac{c(\mathrm{d}x_n)}{4i(\xi_n-i)}+\frac{c(\mathrm{d}x_n)-ic(\xi')}{8(\xi_n-i)^2}+\frac{3\xi_n-7i}{8(\xi_n-i)^3}[ic(\xi')-c(\mathrm{d}x_n)]\right];\nonumber\\
P_1&=\sigma_0(D)\left[\frac{[c(\xi')+ic(\mathrm{d}x_n)]c(\mathrm{d}x_n)}{4(\xi_n-i)}
  +\frac{[c(\xi')+ic(\mathrm{d}x_n)]^2}{4(\xi_n-i)^2}\right];\nonumber\\
  P_2&=c(\mathrm{d}x_n)\partial_{x_n}[c(\xi')]\left[\frac{i[c(\xi')+ic(\mathrm{d}x_n)]c(\mathrm{d}x_n)}{-4(\xi_n-i)}
  +\frac{[c(\xi')+ic(\mathrm{d}x_n)]}{4(\xi_n-i)^2}\right] \nonumber\\
  &+h'(0)|\xi'|^2\bigg[\frac{i[c(\xi')+ic(\mathrm{d}x_n)]^2-2c(\mathrm{d}x_n)[3c(\xi')+4ic(\mathrm{d}x_n)]}{16(\xi_n-i)}\nonumber\\ &+\frac{[c(\xi')+ic(\mathrm{d}x_n)][c(\xi')+5ic(\mathrm{d}x_n)]}{16(\xi_n-i)^2}-\frac{i[c(\xi')+ic(\mathrm{d}x_n)]^2}{8(\xi_n-i)^3}\bigg].\nonumber\\ 
\end{align}
We note that $i<n,~\int_{|\xi'|=1}\xi_{i_{1}}\xi_{i_{2}}\cdots\xi_{i_{2d+1}}\sigma(\xi')=0$,
so we omit some items that have no contribution for computing {\bf case b)}.\\
Then, we have
\begin{align}\label{c22}
\widetilde{\Phi}_4&=Vol(S^{n-2})2^{\frac{n}{2}}\bigg\{-X(Y_n)\frac{2\pi i}{(\frac{n}{2}+2)!}G_0+\frac{1}{n-1}g(X^T,Y^T)h'(0)\frac{\pi i}{4(\frac{n}{2}+3)!}G_1\nonumber\\
&+X_nY_nh'(0)\frac{\pi i}{4(\frac{n}{2}+2)!}G_2\bigg\}dx',\nonumber\\
\end{align} 
where
\begin{align}
G_0&=\left[\frac{3n+2)\xi_n^2+ni\xi_n+2}{(\xi_n+i)^{(\frac{n}{2}+1)}}\right]
 ^{(\frac{n}{2}+2)}\bigg|_{\xi_n=i}\nonumber\\
 &=i^{-n-5}2^{-n-1}\bigg((3n+2)C_{-\frac{n}{2}-1}^{\frac{n}{2}}+(\frac{7n}{2}+2)C_{-\frac{n}{2}-1}^{\frac{n}{2}+1}+nC_{-\frac{n}{2}-1}^{\frac{n}{2}+2}\bigg)(\frac{n}{2}+2)!;\nonumber\\
G_1&=\left[\frac{(1-n)i\xi_n^4+(7-5n)\xi_n^3+(7n-3)i\xi_n^2+(3n+10)\xi_n-8i}{(\xi_n+i)^{\frac{n}{2}+1}}\right]
 ^{(\frac{n}{2}+3)}\bigg|_{\xi_n=i}\nonumber\\
 &=i^{-n-5}4^{-\frac{n}{2}-2}\bigg(16(n-1)C_{-\frac{n}{2}-1}^{\frac{n}{2}-1}+(8n+24)C_{-\frac{n}{2}-1}^{\frac{n}{2}}+(48-8n)C_{-\frac{n}{2}-1}^{\frac{n}{2}+1}+2C_{-\frac{n}{2}-1}^{\frac{n}{2}+2}+C_{-\frac{n}{2}-1}^{\frac{n}{2}+3}\bigg)(\frac{n}{2}+3)!;\nonumber\\
 G_2&=\left[\frac{2n\xi_n^4+(4ni-8i-2)\xi_n^3+(2i+1-5n)\xi_n^2-(3ni-12n+12)\xi_n+13}{(\xi_n+i)^{(\frac{n}{2}+1)}}\right]
 ^{(\frac{n}{2}+2)}\bigg|_{\xi_n=i}\nonumber\\
 &=i^{-n-3}2^{-n-1}\bigg(8nC_{-\frac{n}{2}-1}^{\frac{n}{2}-2}+(24n-16+4i)C_{-\frac{n}{2}-1}^{\frac{n}{2}-1}+(29n-25+4i)C_{-\frac{n}{2}-1}^{\frac{n}{2}}\nonumber\\
 &+(\frac{7+6i}{2}n-3i+1)C_{-\frac{n}{2}-1}^{\frac{n}{2}+2}\bigg)(\frac{n}{2}+2)!.\nonumber\\
\end{align}
\noindent {\bf  case c)}~$r=1,~\ell=-n,~k=j=|\alpha|=0$.\\
By (\ref{c2}), we get
\begin{align}\label{c23}
\widetilde{\Phi}_5&=-i\int_{|\xi'|=1}\int^{+\infty}_{-\infty}{\rm trace} [\pi^+_{\xi_n}\sigma_{1}(\nabla_X^{S(TM)}\nabla_Y^{S(TM)}D^{-1})\times
\partial_{\xi_n}\sigma_{-n}(D^{-(n-1)})](x_0)d\xi_n\sigma(\xi')dx'\nonumber\\
&=i\int_{|\xi'|=1}\int^{+\infty}_{-\infty}{\rm trace} [\partial_{\xi_n}\pi^+_{\xi_n}\sigma_{1}(\nabla_X^{S(TM)}\nabla_Y^{S(TM)}D^{-1})\times
\sigma_{-n}(D^{-(n-1)})](x_0)d\xi_n\sigma(\xi')dx'.
\end{align}
By (3.37) in \cite{Wa7}, we have
\begin{align}\label{c24}
\sigma_{-n}(D^{-(n-1)})(x_0)
&=\frac{(-n-1)h'(0)c(\mathrm{d}x_n)}{4(1+\xi_n^2)^{\frac{n}{2}}}-n\xi_n(1+\xi_n^2)^{-\frac{n}{2}-1}\partial_{x_n}(c(\xi')(x_0)+\frac{n}{2}i(1+\xi_n^2)^{-\frac{n}{2}+1}\nonumber\\
&[c(\xi')+\xi_nc(\mathrm{d}x_n)]\times\bigg[\frac{-ih'(0)c(\xi')c(\mathrm{d}x_n)-(n+1)ih'(0)c(\xi')}{2(1+\xi_n^2)^2}-\frac{2ih'(0)\xi_n}{(1+\xi_n^2)^3}\bigg]\nonumber\\
&-[c(\xi')+\xi_nc(\mathrm{d}x_n)]h'(0)\xi_n[(-\frac{n}{2})^2+\frac{n}{2}][(1+\xi_n^2)^{-\frac{n}{2}-2}].\nonumber\\
\end{align}
By Lemma \ref{lem9356}, we have
\begin{align}\label{c25}
\partial_{\xi_n}\pi^+_{\xi_n}\sigma_{1}(\nabla_X^{S(TM)}\nabla_Y^{S(TM)}D^{-1})
&=\frac{c(\xi')+ic(\mathrm{d}x_n)}{2(\xi_n-i)^2}\bigg(\sum_{j,l=1}^{n-1}X_jY_l\xi_j\xi_l
-X_nY_n\bigg) \nonumber\\
&+\frac{ic(\xi')-c(\mathrm{d}x_n)}{2(\xi_n-i)^2}\sum_{j=1}^{n-1}(X_jY_n
+X_nY_j)\xi_j.
\end{align}
We note that $i<n,~\int_{|\xi'|=1}\xi_{i_{1}}\xi_{i_{2}}\cdots\xi_{i_{2d+1}}\sigma(\xi')=0$,
so we omit some items that have no contribution for computing {\bf case c)}.\\
Therefore, we have
\begin{align}\label{c27}
\widetilde{\Phi}_5&=Vol(S^{n-2})\bigg(X_nY_n-\frac{1}{n-1}g(X^T,Y^T)\bigg)h'(0)\bigg((n+1)2^{\frac{n}{2}-3}\frac{2\pi i}{(\frac{n}{2}+1)!}H_0-2^{\frac{n}{2}-1}\frac{2n\pi}{(\frac{n}{2}+2)!}H_1\nonumber\\
&-2^{\frac{n}{2}-2}\frac{2\pi}{(\frac{n}{2}+2)!}H_2-2^{\frac{n}{2}}\frac{n\pi}{(\frac{n}{2}+3)!}H_3\bigg)dx',
\end{align}
where
\begin{align}
H_0&=\left[\frac{1}{(\xi_n+i)^{\frac{n}{2}}}\right]
 ^{(\frac{n}{2}+1)}\bigg|_{\xi_n=i}=(2i)^{-n-1}A_{-\frac{n}{2}}^{\frac{n}{2}+1};\nonumber\\
H_1&=\left[\frac{\xi_n}{(\xi_n+i)^{(\frac{n}{2}+1)}}\right]
 ^{(\frac{n}{2}+2)}\bigg|_{\xi_n=i}=(2i)^{-n-3}\bigg(2iC_{-\frac{n}{2}-1}^{\frac{n}{2}+1}+iC_{-\frac{n}{2}-1}^{\frac{n}{2}+2}\bigg)(\frac{n}{2}+2)!;\nonumber\\
H_2&=\left[\frac{-(n+1)ni\xi_n^2-n^2\xi_n+ni}{(\xi_n+i)^{(\frac{n}{2}+1)}}\right]
 ^{(\frac{n}{2}+2)}\bigg|_{\xi_n=i}\nonumber\\
 &=-i^{-n-4}4^{-\frac{n}{2}-1}\bigg((4n-4i+4)C_{-\frac{n}{2}-1}^{\frac{n}{2}}+(4n-4i+4)C_{-\frac{n}{2}-1}^{\frac{n}{2}+1}+(n-i+2)C_{-\frac{n}{2}-1}^{\frac{n}{2}+2}\bigg)(\frac{n}{2}+2)!;\nonumber\\
H_3&=\left[\frac{i\xi_n^2+\xi_n}{(\xi_n+i)^{(\frac{n}{2}+2)}}\right]
 ^{(\frac{n}{2}+3)}\bigg|_{\xi_n=i}=i^{-n-6}4^{-\frac{n}{2}-2}\bigg(2C_{-\frac{n}{2}-2}^{\frac{n}{2}+1}+C_{-\frac{n}{2}-2}^{\frac{n}{2}+1}\bigg)(\frac{n}{2}+3)!.\nonumber\\
\end{align}
 Now $\widetilde{\Phi}$ is the sum of the cases (a), (b) and (c). Therefore, we get
\begin{align}\label{c28}
\widetilde{\Phi}&=\sum_{i=1}^5\widetilde{\Phi}_i\nonumber\\
&=Vol(S^{n-2})2^{\frac{n}{2}}\bigg\{-\frac{1}{n-1}\bigg(h'(0)g(X^T,Y^T)+\partial_{x_n}(g(X^T,Y^T))\bigg)\frac{\pi i}{(\frac{n}{2}+2)!}E_{0}-\partial_{x_n}(X_nY_n) \frac{\pi i}{(\frac{n}{2}+2)!}E_{1}\nonumber\\
&-(i-1)h'(0)\bigg(\frac{1}{n-1}g(X^T,Y^T)-X_nY_n\bigg)\frac{\pi i}{2(\frac{n}{2}+2)!}E_{2}-\frac{1}{n-1}g(X^T,Y^T)h'(0)\frac{\pi i}{2(\frac{n}{2}+3)!}E_{3}+X_nY_n\nonumber\\
&\frac{\pi i}{(\frac{n}{2}+3)!}E_{4}-\bigg(\frac{1}{n-1}g(X^T,Y^T)+X_nY_n\bigg)\frac{2\pi i}{(\frac{n}{2}+3)!}F_0-X(Y_n)\frac{2\pi i}{(\frac{n}{2}+2)!}G_0+\frac{1}{n-1}g(X^T,Y^T)h'(0)\nonumber\\
&\frac{\pi i}{4(\frac{n}{2}+3)!}G_1+X_nY_nh'(0)\frac{\pi i}{4(\frac{n}{2}+2)!}G_2+h'(0)\bigg(X_nY_n-\frac{1}{n-1}g(X^T,Y^T)\bigg)\bigg(\frac{n+1}{3}\frac{\pi i}{(\frac{n}{2}+1)!}H_0-\frac{n\pi}{(\frac{n}{2}+2)!}\nonumber\\
&H_1-\frac{\pi }{2(\frac{n}{2}+2)!}H_2-\frac{2n\pi}{(\frac{n}{2}+3)!}H_3\bigg)\bigg\}dx'.\nonumber\\
\end{align}
Then, by (\ref{c28}), we obtain following theorem
\begin{thm}\label{cthmb1}
Let $M$ be an $n$-dimensional oriented
compact spin manifold with boundary $\partial M$ and n is even, then we get the following equality:
\begin{align}
\label{cb2773}
&\widetilde{{\rm Wres}}[\pi^+(\nabla_X^{S(TM)}\nabla_Y^{S(TM)}D^{-1})\circ\pi^+(D^{-(n-1)})]\nonumber\\
&=\frac{(2\pi)^{\frac{n}{2}}}{3(\frac{n}{2}-1)!}\int_{M}[Ric(X,Y)-\frac{1}{2}sg(X,Y)]d{\rm Vol_{M}}+\frac{(2\pi)^{\frac{n}{2}}}{4(\frac{n}{2}-1)!}\int_{M}sg(X,Y)d{\rm Vol_{M}}+\int_{\partial M}Vol(S^{n-2})2^{\frac{n}{2}}\nonumber\\
&\bigg\{-\frac{1}{n-1}\bigg(h'(0)g(X^T,Y^T)+\partial_{x_n}(g(X^T,Y^T))\bigg)\frac{\pi i}{(\frac{n}{2}+2)!}E_{0}-\partial_{x_n}(X_nY_n) \frac{\pi i}{(\frac{n}{2}+2)!}E_{1}\nonumber\\
&-(i-1)h'(0)\bigg(\frac{1}{n-1}g(X^T,Y^T)-X_nY_n\bigg)\frac{\pi i}{2(\frac{n}{2}+2)!}E_{2}-\frac{1}{n-1}g(X^T,Y^T)h'(0)\frac{\pi i}{2(\frac{n}{2}+3)!}E_{3}+X_nY_n\nonumber\\
&\frac{\pi i}{(\frac{n}{2}+3)!}E_{4}-\bigg(\frac{1}{n-1}g(X^T,Y^T)+X_nY_n\bigg)\frac{2\pi i}{(\frac{n}{2}+3)!}F_0-X(Y_n)\frac{2\pi i}{(\frac{n}{2}+2)!}G_0+\frac{1}{n-1}g(X^T,Y^T)h'(0)\nonumber\\
&\frac{\pi i}{4(\frac{n}{2}+3)!}G_1+X_nY_nh'(0)\frac{\pi i}{4(\frac{n}{2}+2)!}G_2+h'(0)\bigg(X_nY_n-\frac{1}{n-1}g(X^T,Y^T)\bigg)\bigg(\frac{n+1}{3}\frac{\pi i}{(\frac{n}{2}+1)!}H_0-\frac{n\pi}{(\frac{n}{2}+2)!}\nonumber\\
&H_1-\frac{\pi }{2(\frac{n}{2}+2)!}H_2-\frac{2n\pi}{(\frac{n}{2}+3)!}H_3\bigg)\bigg\}d{\rm Vol_{M}}.
\end{align}
\end{thm}
\section{Dabrowski-Sitarz-Zalecki type theorems for odd dimensional manifolds with boundary}
\label{section:4}
For an $n-1$ dimensional manifolds with boundary and $n-1$ is odd, as in section 5,6 and 7 in \cite{Wa1}, we have the formula\\
\begin{align}\label{b1}
\widetilde{{\rm Wres}}[\pi^+(\nabla_X^{S(TM)}\nabla_Y^{S(TM)}D^{-2})\circ\pi^+(D^{-(n-2)})]=\int_{\partial M}\Psi.
\end{align}
Then in (\ref{qqq15}), $r-k-|\alpha|+l-j-1=-(n-1)$, $r,l \leq 0,$ so we get $r=0,~l=-(n-2),k=|\alpha|=j=0,$ then
\begin{eqnarray}\label{b2}
\Psi&=\int_{|\xi'|=1}\int^{+\infty}_{-\infty}{\rm trace}_{S(TM)}[\pi^+_{\xi_n}\sigma_{0}(\nabla_X^{S(TM)}\nabla_Y^{S(TM)}D^{-2})
\nonumber\\
&\times\partial_{\xi_n}\sigma_{-(n-2)}(D^{-(n-2)})(x',0,\xi',\xi_n)]d\xi_3\sigma(\xi')dx'.
\end{eqnarray}
By (\ref{aa38}) and (\ref{62})and the Cauchy integral formula, we get\\
\begin{align}\label{b5}
\Psi&=-Vol(S^{n-2})(\frac{\sum_{j=1}^{n-1}(X_jY_j}{n-1}+X_nY_n)(1-\frac{n}{2})2^{\frac{n}{2}-1}\frac{2\pi i}{\frac{n}{2}!}\left[\frac{\xi_{n}}{(\xi_n+i)^{\frac{n}{2}}}\right]^{(\frac{n}{2})}\bigg|_{\xi_n=i}dx'\nonumber\\
&:=-Vol(S^{n-2})\bigg(\frac{1}{n-1}g(X^T,Y^T)+X_nY_n\bigg)(1-\frac{n}{2})2^{\frac{n}{2}-1}\frac{2\pi i}{\frac{n}{2}!}M_0dx',
\end{align}
where
\begin{align}
M_0&=\left[\frac{\xi_{n}}{(\xi_n+i)^{\frac{n}{2}}}\right]^{\frac{n}{2}}\bigg|_{\xi_n=i}=i^{-n+1}2^{-n}\bigg(2C_{-\frac{n}{2}}^{\frac{n}{2}-1}+C_{-\frac{n}{2}}^{\frac{n}{2}}\bigg)\frac{n}{2}!.\nonumber\\
\end{align}
Therefore, we get the following theorem
\begin{thm}\label{athm3}
Let $M$ be an $n-1$ dimensional
oriented compact spin manifold with boundary $\partial M$ and $n-1$ is odd, then
\begin{align}\label{b7}
&\widetilde{{\rm Wres}}[\pi^+(\nabla_X^{S(TM)}\nabla_Y^{S(TM)}D^{-2})\circ\pi^+(D^{-(n-2)})]\nonumber\\
&=\int_{\partial M}-Vol(S^{n-2})\bigg(\frac{1}{n-1}g(X^T,Y^T)+X_nY_n\bigg)(1-\frac{n}{2})2^{\frac{n}{2}-1}\frac{2\pi i}{\frac{n}{2}!}M_0d{\rm Vol_{M}}\nonumber\\
\end{align}
where $vol_{\partial M}$ denotes the canonical volume form of $\partial M$.\\
\end{thm}
For an $n-1$ dimensional manifolds with boundary, as in section 5,6 and 7 in \cite{Wa1}, we have the formula\\
\begin{align}\label{b1}
\widetilde{{\rm Wres}}[\pi^+(\nabla_X^{S(TM)}\nabla_Y^{S(TM)}D^{-1})\circ\pi^+(D^{-(n-1)})]=\int_{\partial M}\widetilde{\Psi}.
\end{align}
Then in (\ref{c2}), $r-k-|\alpha|+l-j-1=-(n-1)$, $r,l \leq 0,$ so we get $r=1,~l=-(n-1),k=|\alpha|=j=0,$ then
\begin{eqnarray}\label{b2}
\widetilde{\Psi}&=\int_{|\xi'|=1}\int^{+\infty}_{-\infty}{\rm trace}_{S(TM)}[\pi^+_{\xi_n}\sigma_{1}(\nabla_X^{S(TM)}\nabla_Y^{S(TM)}D^{-1})
\nonumber\\
&\times\partial_{\xi_n}\sigma_{-(n-1)}(D^{-(n-1)})(x',0,\xi',\xi_n)]d\xi_3\sigma(\xi')dx'.
\end{eqnarray}
Moreover,
\begin{align}\label{c38}
\pi^+_{\xi_n}\sigma_{1}(\nabla_X^{S(TM)}\nabla_Y^{S(TM)}D^{-1})
&=-\frac{c(\xi')+ic(\mathrm{d}x_n)}{2(\xi_n-i)}(\sum_{j,l=1}^{n-1}X_jY_l\xi_j\xi_l
+X_nY_n)\nonumber\\
&-\frac{ic(\xi')-c(\mathrm{d}x_n)}{2(\xi_n-i)}\sum_{j=1}^{n-1}(X_jY_n
+X_nY_j)\xi_j.
\end{align}

By (\ref{c20}) and (\ref{c38}) and the Cauchy integral formula, we get\\
\begin{align}\label{b5}
\widetilde{\Psi}&=-\frac{i}{2}h'(0)Vol(S^{n-2})\bigg(\frac{\sum_{j=1}^{n-1}X_jY_j}{n-1}+X_nY_n\bigg)2^{\frac{n}{2}+1}\frac{2\pi i}{(\frac{n}{2}+1)!}\left[\frac{ni\xi_{n}+ni}{(\xi_n+i)^{(\frac{n}{2}+1)}}\right]^{(\frac{n}{2}+1)}\bigg|_{\xi_n=i}dx'\nonumber\\
&-ih'(0)Vol(S^{n-2})\bigg(\frac{\sum_{j=1}^{n-1}X_jY_j}{n-1}+X_nY_n\bigg)2^{\frac{n}{2}}\frac{2\pi i}{\frac{n}{2}!}\left[\frac{1}{(\xi_n+i)^{\frac{n}{2}}}\right]^{(\frac{n}{2})}\bigg|_{\xi_n=i}dx'\nonumber\\
&:=-iVol(S^{n-2})h'(0)\bigg(\frac{1}{n-1}g(X^T,Y^T)+X_nY_n\bigg)2^{\frac{n}{2}}\bigg(\frac{\pi i}{2(\frac{n}{2}+1)!}N_0+
\frac{\pi i}{\frac{n}{2}!}N_1\bigg)dx',
\end{align}
where
\begin{align}
N_0&=\left[\frac{ni\xi_{n}+ni}{(\xi_n+i)^{(\frac{n}{2}+1)}}\right]^{(\frac{n}{2}+1)}\bigg|_{\xi_n=i}=(2i)^{-n-1}\bigg((2ni+2i)C_{-\frac{n}{2}-1}^{\frac{n}{2}}+(ni+\frac{n}{2}+i)C_{-\frac{n}{2}-1}^{\frac{n}{2}+1}\bigg)(\frac{n}{2}+1)!;\nonumber\\
N_1&=\left[\frac{1}{(\xi_n+i)^{\frac{n}{2}}}\right]^{\frac{n}{2}}\bigg|_{\xi_n=i}=(2i)^{-n}A_{-\frac{n}{2}}^{\frac{n}{2}}.\nonumber\\
\end{align}
Therefore, we get the following theorem
\begin{thm}\label{athm3}
Let $M$ be an $n-1$ dimensional
oriented compact spin manifold with boundary $\partial M$ and $n-1$ is odd, then
\begin{align}\label{b7}
&\widetilde{{\rm Wres}}[\pi^+(\nabla_X^{S(TM)}\nabla_Y^{S(TM)}D^{-1})\circ\pi^+(D^{-(n-1)})]\nonumber\\
&=\int_{\partial M}-iVol(S^{n-2})h'(0)2^{\frac{n}{2}}\bigg(\frac{1}{n-1}g(X^T,Y^T)+X_nY_n\bigg)\bigg(\frac{\pi i}{2(\frac{n}{2}+1)!}N_0+
\frac{\pi i}{\frac{n}{2}!}N_1\bigg)d{\rm Vol_{M}}\nonumber\\
\end{align}
where $vol_{\partial M}$ denotes the canonical volume form of $\partial M$.\\
\end{thm}
\section*{Acknowledgements}
This work was supported by NSFC. 11771070 .
 The authors thank the referee for his (or her) careful reading and helpful comments.

\section*{}

\end{document}